\documentclass[11pt]{article}
\usepackage[T1]{fontenc}
\usepackage{lmodern}
\usepackage[margin=1in]{geometry}
\usepackage{amsmath,amssymb,amsthm,mathtools,mathrsfs}
\usepackage{microtype,enumitem}
\usepackage{xcolor,xurl}
\usepackage{hyperref}
\hypersetup{colorlinks=true,linkcolor=teal!65!black,citecolor=teal!65!black,
urlcolor=teal!65!black,
pdfauthor={Yi Wang and Shuhan Yang},
pdftitle={Higher-Order Minkowski Inequalities via Gradient-Weighted Hessian Potentials},
pdfsubject={Newton flux, conformal divergence identities, and higher-order Minkowski inequalities}}
\numberwithin{equation}{section}
\allowdisplaybreaks[2]
\setlist[enumerate]{itemsep=4pt,topsep=6pt}
\setlist[itemize]{itemsep=3pt,topsep=5pt}
\newtheorem{theorem}{Theorem}[section]
\newtheorem{proposition}[theorem]{Proposition}
\newtheorem{lemma}[theorem]{Lemma}

\theoremstyle{definition}
\newtheorem{definition}[theorem]{Definition}
\newtheorem{assumption}[theorem]{Assumption}
\theoremstyle{remark}
\newtheorem{remark}[theorem]{Remark}
\newtheorem*{remark*}{Remark}
\theoremstyle{plain}
\newtheorem*{introtheorem}{Theorem A}
\newtheorem*{introcorollary}{Theorem B}
\newtheorem*{mainresult}{Theorem A}
\newcommand{\RR}{\mathbb R}
\newcommand{\Sphere}{\mathbb S}
\newcommand{\Div}{\operatorname{div}}
\newcommand{\tr}{\operatorname{tr}}
\newcommand{\Crit}{\operatorname{Crit}}

\newcommand{\LL}{\mathcal L_{k,p}}
\newcommand{\KK}{\mathcal K}
\newcommand{\HH}{\mathcal H}
\newcommand{\CC}{\mathcal C_{k,p}}
\newcommand{\dd}{\,\mathrm d}
\newcommand{\ip}[2]{\langle #1,#2\rangle}

\newcommand{\avgint}{\frac{1}{\omega_{n-1}}\int}
\title{Higher-Order Minkowski Inequalities\\
via Gradient-Weighted Hessian Potentials}
\author{Yi Wang 
\thanks{3400 N Charles St, Baltimore, MD, 21218, Department of Mathematics, Johns Hopkins University, ywang261@jhu.edu}
, Shuhan Yang \thanks{3400 N Charles St, Baltimore, MD, 21218, Department of Mathematics, Johns Hopkins University, syang180@jhu.edu}}
\date{}
\begin{document}
\maketitle
\begin{abstract}
We prove a sharp boundary inequality for exterior potentials of a
gradient-weighted $k$-Hessian equation under regularity, level-curvature,
flux conservation, and radial decay hypotheses. The inequality bounds
a power of the Newton flux by the normalized integral of
$H_k^p/H_{k-1}^{p-1}$, and equality characterizes balls. After a
conformal change of metric, two Newton-tensor divergence
identities give comparisons between regular levels even when the slab
between them contains critical points. For smooth strictly star-shaped
domains with strictly $k$-convex boundary, a companion existence theorem
supplies exterior potentials for $p$ close to one. We identify the
limiting flux with the integral of $H_{k-1}$ and recover the sharp
consecutive Minkowski inequality, including its equality case.
\end{abstract}

\section{Introduction}\label{sec:introduction}

The Minkowski inequalities relate curvature integrals of the boundary of
a domain and characterize Euclidean balls through equality. Let
\(\Omega\subset\RR^n\) be a smooth bounded domain, let \(h_\delta\) denote
the outward shape operator of its boundary, and write
\[
 H_j=\binom{n-1}{j}^{-1}\sigma_j(h_\delta),\qquad
 A_j(\Omega)=\frac{1}{\omega_{n-1}}
       \int_{\partial\Omega}H_j\dd S_\delta,
 \qquad \omega_{n-1}=|\Sphere^{n-1}|.
\]
The consecutive curvature-integral inequality takes the form
\[
 A_k(\Omega)\ge
 A_{k-1}(\Omega)^{\frac{n-k-1}{n-k}}.
\]
For convex domains this is an Alexandrov--Fenchel inequality for mixed
volumes \cite{Ale37,Ale38,hormander1994notions}. Trudinger studied
isoperimetric inequalities for quermassintegrals beyond the convex case
\cite{TRUDINGER1994411}. Guan and Li~\cite{GuanLi2009} established the
sharp inequality, together with its equality characterization, for
smooth star-shaped domains with \(k\)-convex boundary by an expanding
curvature flow of the type studied by Gerhardt and Urbas
\cite{Gerhardt1990,Urbas1990}. Level-set formulations were developed
by Evans and Spruck for mean curvature flow \cite{Evans1991635} and
by Huisken and Ilmanen for inverse mean curvature flow \cite{HuiskenIlmanen2001}.

The classical isoperimetric inequality admits both geometric proofs
\cite{Chavel} and an elliptic proof via the Alexandrov--Bakelman--Pucci
estimate \cite{Cabre1,Cabre2}. On submanifolds, the Michael--Simon
Sobolev inequality involves the mean curvature \cite{Michael1973SobolevAM}.
Castillon obtained related inequalities by optimal transport
\cite{CASTILLON201079}, and Wang proved analogues for higher mean
curvatures \cite{wang2013michaelsimoninequalitieskthmean}. In codimension
at most two, Brendle proved the sharp submanifold Sobolev inequality
\cite{brendle2020isoperimetricinequalityminimalsubmanifold}, for which
Brendle and Eichmair later gave a transport proof
\cite{brendle2023proofmichaelsimonsobolevinequalityusing}.

Chang and Wang used optimal transport to study
quermassintegrals on \(k\)-convex domains \cite{ChangWang2013}.
They also obtained bounds on volume in terms of curvature integrals
\cite{ChangWang2014}, which Qiu extended to higher orders \cite{Qiu2015}.
For nearly spherical domains, Fuglede proved stability of the classical
isoperimetric inequality \cite{Fuglede1989}. Related Minkowski inequalities are proved
in \cite{Glaudo2022,WangYang2025} and quantitative quermassintegral
inequalities in \cite{VanBlarganWangQuant2024,VanBlarganWangFlow2024}.
Another approach to stability uses almost-umbilicity estimates
\cite{De_Rosa_2021,Scheuer+2025+381+405}.

A different approach to Minkowski inequalities uses the level
hypersurfaces of exterior potentials. Monotonicity formulas along these
levels were obtained for harmonic potentials in \cite{AgostinianiMazzieri2020}
and for \(p\)-capacitary potentials of convex domains in \cite{FMP2019}.
Agostiniani, Fogagnolo, and
Mazzieri~\cite{AFM} developed effective monotonicity formulas for
\(p\)-capacitary potentials and used the limit \(p\downarrow1\) to derive
Minkowski-type inequalities. Their method accommodates critical points
of the potential through divergence identities and cutoff arguments.
The purpose of the present paper is to develop a higher-order version
of this method for the gradient-weighted Hessian operator
\[
 \LL u=\Div_\delta\!\left(
 |Du|^{p-k-1}T_{k-1}(D^2u)Du\right),
\]
where \(T_{k-1}\) is the \((k-1)\)-st Newton tensor. This operator is the
\(p\)-Laplacian when \(k=1\), and is \(k\sigma_k(D^2u)\) when \(p=k+1\).
Tso related Hessian equations to curvature integrals through
symmetrization \cite{TSO89}. For the homogeneous \(k\)-Hessian equation, Ma and Zhang studied the
exterior Dirichlet problem \cite{MaZhang2024}, and Xiao derived a
Minkowski inequality from its exterior level sets \cite{Xiao2022}.

Fix
\[
 n\ge3,\qquad 1\le k\le n-2,\qquad 1<p<n-k,\qquad
 \alpha_p=\frac{n-k-p+1}{p-1}.
\]
Write $E=\RR^n\setminus\overline\Omega$. The normalized boundary flux
of an exterior potential $u$ is
\[
 \CC(\Omega;u)=\frac{1}{\omega_{n-1}}
       \int_{\partial\Omega}
       \left(\frac{|Du|}{\alpha_p}\right)^{p-1}
       H_{k-1}\dd S_\delta.
\]

\noindent\hypertarget{intro-potential-hypotheses}{\textbf{Hypotheses on the exterior potential in Theorem~A.}}
\begin{enumerate}[label=\textup{(P\arabic*)},leftmargin=*,itemsep=2pt]
\item $u\in C^1_{\mathrm{loc}}(\overline E)$, $0<u<1$ in $E$,
$u=1$ on $\partial\Omega$, and $u(x)\to0$ as $|x|\to\infty$.
\item On $\{x\in E:Du(x)\ne0\}$, $u$ is $C^3$ and solves $\LL u=0$
classically. It is also $C^3$ up to $\partial\Omega$ in a boundary
collar, with $|Du|>0$ on $\partial\Omega$.
\item At each point where $Du\ne0$ in $E$, the level hypersurface
oriented by $-Du/|Du|$ has $\sigma_{k-1}(h_\delta)>0$, where
$h_\delta$ is its Euclidean shape operator.
\item For $z=-u$, the field $|Dz|^{p-k-1}T_{k-1}(D^2z)Dz$ on
$\{Dz\ne0\}$ admits an $L^1_{\mathrm{loc}}(E)$ extension whose
distributional divergence vanishes throughout $E$.
\item For some $\mathsf a>0$, uniformly in direction as $|x|\to\infty$,
\[
 D^j\bigl(u(x)-\mathsf a|x|^{-\alpha_p}\bigr)
 =o(|x|^{-\alpha_p-j}),\qquad j=0,1,2.
\]
\end{enumerate}

\begin{introtheorem}
Let \(\Omega\subset\RR^n\) be a smooth bounded strictly star-shaped
domain with strictly \(k\)-convex boundary. If $u$ satisfies
\hyperlink{intro-potential-hypotheses}{\textup{(P1)--(P5)}} above, then
\[
 \CC(\Omega;u)^{\frac{n-k-p}{n-k-p+1}}
 \le \frac{1}{\omega_{n-1}}\int_{\partial\Omega}
       \frac{H_k^p}{H_{k-1}^{p-1}}\dd S_\delta.
\]
Equality holds if and only if $\Omega$ is a ball. If
$\Omega=B_R(x_0)$, the potential is necessarily
$u(x)=(R/|x-x_0|)^{\alpha_p}$.
\end{introtheorem}

For each fixed domain $\Omega$ as in Theorem~A,
\cite[Theorem~1.1]{ExteriorConstruction} constructs a smooth exterior
potential $u_p$ satisfying all five hypotheses \textup{(P1)--(P5)}
for every $1<p<1+\varepsilon_0(\Omega,n,k)$. Thus Theorem~A applies to
each $u_p$. The same construction gives the uniform boundary limit
$(p-1)|Du_p|\to(n-k)H_k/H_{k-1}$, which identifies the limiting flux
and yields the curvature inequality in Theorem~B upon letting $p\downarrow1$.

\begin{introcorollary}
Let $n\ge3$, $1\le k\le n-2$, and let $\Omega\subset\RR^n$ be a
smooth bounded strictly star-shaped domain with strictly $k$-convex
boundary. For the exterior potentials $u_p$ constructed in
\cite{ExteriorConstruction},
\[
 \lim_{p\downarrow1}\CC(\Omega;u_p)=A_{k-1}(\Omega),\qquad
 A_k(\Omega)\ge A_{k-1}(\Omega)^{(n-k-1)/(n-k)}.
\]
Equality in the curvature inequality holds if and only if $\Omega$ is a ball.
\end{introcorollary}

The proof of Theorem~A takes place in the conformal metric associated
with the logarithm of $u$. Two divergence identities reduce the sign
of the relevant level-set fluxes to scalar and tensorial Newton
inequalities. A cutoff in the gradient variable leaves the flux on both
regular endpoints intact and vanishes near every intervening critical
point. Distributional conservation of the original Euclidean flux and
the radial asymptotics then identify its value at infinity.
At equality the Newton deficit vanishes on the boundary, forcing it to
be umbilic.

Section~\ref{sec:notation} fixes the curvature conventions.
Section~\ref{sec:preparatory} identifies the conserved flux and establishes
the conformal identities. Section~\ref{sec:proof} uses these identities
to prove Theorem~A. Section~\ref{sec:limit} deduces Theorem~B from
Theorem~A using the interior-ball comparison of \cite{ExteriorConstruction}
and a radial barrier.

Appendix~\ref{app:ellipticity} computes the principal linearization,
and Appendix~\ref{app:energy-detailed} gives its variational counterpart
for the exterior equation.

\section{Notation and geometric assumptions}\label{sec:notation}
Fix
\begin{equation*}
 n\ge3,\qquad 1\le k\le n-2,\qquad 1<p<n-k.
\end{equation*}
Thus \(n-k\ge2\) and \(n-k-p+1>1\). The ambient space is \(\RR^n\),
and every regular level hypersurface has dimension \(n-1\).
Let \(\Omega\) be a bounded smooth strictly star-shaped domain with
strictly \(k\)-convex boundary, and put
\(E=\RR^n\setminus\overline\Omega\). After translating a star center
to the origin, the geometric assumptions mean
\[
 \partial\Omega=\{r_\Omega(\theta)\theta:\theta\in\Sphere^{n-1}\},
 \qquad r_\Omega\in C^\infty(\Sphere^{n-1}),\quad r_\Omega>0,
\]
\[
 \sigma_j(\kappa)>0\quad\hbox{on }\partial\Omega,\qquad 1\le j\le k,
\]
where \(\kappa\) denotes the outward principal curvature vector.
In particular \(x\cdot\nu>0\) on the boundary.

For a self-adjoint endomorphism \(A\), \(\sigma_j(A)\) is the
\emph{unnormalized} elementary symmetric polynomial. For vectors in
\(\RR^m\), with \(1\le j\le m\), write
\[
 \Gamma_j=\{\lambda:\sigma_i(\lambda)>0\text{ for }1\le i\le j\};
\]
membership of a symmetric matrix means membership of its eigenvalue
vector. We use
\[
 T_0(A)=I,\qquad T_j(A)=\sigma_j(A)I-A T_{j-1}(A),
 \qquad T_{-1}=0,\qquad \sigma_0=1.
\]
For an $m\times m$ matrix $A$, we set $\sigma_j(A)=0$ for $j>m$.
Newton tensors of ambient matrices act on $\RR^n$, whereas those of
hypersurface shape operators act on the $(n-1)$-dimensional tangent spaces.
Euclidean derivatives and the Euclidean metric are denoted by \(D\)
and \(\delta\). The operator under consideration is
\begin{equation*}
 \LL w=\Div_\delta\!\left(
 |Dw|_\delta^{p-k-1}T_{k-1}(D_\delta^2w)Dw\right).
\end{equation*}
It equals \(\Delta_p\) for \(k=1\). The ordinary \(k\)-Hessian
operator, multiplied by \(k\), occurs at \(p=k+1\).
The outward shape operator is \(h(V)=\nabla_V\nu\), so the principal
curvatures of a Euclidean sphere of radius \(R\) are \(1/R\).
Set
\[
 \omega_{n-1}=|\Sphere^{n-1}|,\qquad
 H_j=\frac{\sigma_j(h_\delta)}{\binom{n-1}{j}},\qquad
 A_j(\Omega)=\avgint_{\partial\Omega}H_j\dd S_\delta.
\]
Here \(H_0=1\) and \(A_j(B_R)=R^{n-1-j}\).

\section{Exterior potentials and conformal identities}\label{sec:preparatory}

\subsection{Exterior potentials and normalized flux}\label{sec:potentials}

With the conventions of Section~\ref{sec:notation}, the exterior problem is
\begin{equation}\label{eq:exterior}
 \LL u=0\quad\text{in }E,\qquad
 u=1\quad\text{on }\partial\Omega,\qquad
 \lim_{|x|\to\infty}u(x)=0.
\end{equation}

\begin{definition}[Normalized Newton flux]\label{def:flux}
Let $u\in C^1_{\mathrm{loc}}(\overline E)$.
The normalized Newton flux of $u$ is defined by
\begin{equation}\label{eq:capacitydef}
 \CC(\Omega;u)
 =\frac{1}{\omega_{n-1}}
 \int_{\partial\Omega}\left(\frac{(p-1)|Du|}{n-k-p+1}\right)^{p-1}
 H_{k-1}\dd S_\delta.
\end{equation}
For a fixed $u$, we write $\CC(\Omega)=\CC(\Omega;u)$.
The flux is attached to the pair $(\Omega,u)$.
\end{definition}

The analytic hypotheses for the exterior potential are stated below.
They extend the regularity and asymptotic properties used in the
$p$-capacitary argument of \cite{AFM} by including level-curvature
positivity and distributional conservation of the Newton flux.
In Section~\ref{sec:limit} the companion construction provides this
input for every sufficiently small $p-1>0$.

\begin{assumption}[Regularity, level curvature, and flux conservation]\label{ass:potential}
For the given domain and parameter, a function $u$ is available with
the following properties.
\begin{enumerate}[label=\textup{(P\arabic*)},ref=P\arabic*]
\item\label{hyp:potential}
$u\in C^1_{\mathrm{loc}}(\overline E)$, $0<u<1$ in $E$, and the
boundary and infinity conditions in \eqref{eq:exterior} hold.
The notation means local regularity up to the finite boundary.
\item\label{hyp:regularity}
On $\mathcal R=\{x\in E:Du(x)\ne0\}$, $u$ is $C^3$ and solves
$\LL u=0$ classically. It is also $C^3$ up to $\partial\Omega$
in a collar, and $|Du|>0$ on $\partial\Omega$.
\item\label{hyp:curvature}
At every point of $\mathcal R$, the local level hypersurface,
oriented by $-Du/|Du|$, satisfies
$\sigma_{k-1}(h_\delta)>0$. This includes noncritical points whose
level value is attained at a critical point elsewhere.
At a globally regular value $\tau\in(0,1)$, the level is oriented as
the boundary of
$\overline\Omega\cup\{x\in E:u(x)\ge\tau\}$.
Full convexity of the levels is not assumed.
\item\label{hyp:flux}
For $z=-u$, the classical vector field
\[
 F_z=|Dz|^{p-k-1}T_{k-1}(D^2z)Dz\qquad\text{on }\mathcal R
\]
admits an $L^1_{\mathrm{loc}}(E)$ extension, still denoted by $F_z$,
with $\Div_\delta F_z=0$ in distributions on $E$.
The extension agrees with the classical field wherever $Dz\ne0$.
\end{enumerate}
\end{assumption}

\begin{assumption}[Radial asymptotics]\label{ass:asymptotic}
\textup{(P5)} There exists $\mathsf a>0$ such that, uniformly in
direction,
\begin{equation}\label{eq:asymptotic}
 D^j\bigl(u(x)-\mathsf a|x|^{-\frac{n-k-p+1}{p-1}}\bigr)
 =o(|x|^{-\frac{n-k-p+1}{p-1}-j}),\qquad j=0,1,2.
\end{equation}
\end{assumption}

Condition \textup{(P5)} implies that all sufficiently small positive
$u$-levels are noncritical. Together with \textup{(P2)}, which excludes
critical points near the boundary, this places $\Crit(u)$ in a compact
subset of $E$. The expansion through second derivatives supplies the
level-set curvature asymptotics used below. It is required separately
at each fixed $p$, with no uniformity as $p\downarrow1$.

The regularity in \textup{(P1)--(P2)} permits the cutoff constructions,
while \textup{(P3)} gives the signs in the geometric comparisons.
These comparisons are stated at globally regular endpoint levels;
their proof requires neither a dimension bound for the critical set
nor a Sard theorem for all critical values. Conditions \textup{(P4)}
and \textup{(P5)} together determine the boundary flux from the radial
asymptotic coefficient. Appendix~\ref{app:ellipticity} discusses the
relation between strict ellipticity and \textup{(P3)}.

\begin{lemma}[Classical flux identity and its normalization]\label{lem:flux}
Under Assumptions \ref{ass:potential}--\ref{ass:asymptotic}, every
globally regular level $\{u=\tau\}$ has outward flux
\begin{equation}\label{eq:flatflux}
 \left\langle F_z,-\frac{Du}{|Du|}\right\rangle_\delta
 =|Du|^{p-1}\sigma_{k-1}(h_\delta).
\end{equation}
Its integral equals the corresponding boundary integral and the
flux through every sufficiently large Euclidean sphere. Moreover,
\begin{equation}\label{eq:Ccoefficient}
 \CC(\Omega)=\mathsf a^{\,p-1}.
\end{equation}
\end{lemma}
\begin{proof}
At a regular point choose a Euclidean orthonormal frame with last
vector $\nu_\delta=Dz/|Dz|$. Differentiating $Dz=|Du|\nu_\delta$
along the level gives $D^2z|_{T\Sigma}=|Du|h_\delta$. In the
principal-minor expansion of $\sigma_k(D^2z)$, differentiation
with respect to the last diagonal entry yields
\[
 T_{k-1}(D^2z)_{nn}=\sigma_{k-1}(D^2z|_{T\Sigma})
                  =|Du|^{k-1}\sigma_{k-1}(h_\delta).
\]
Contracting $F_z$ with $\nu_\delta$ proves \eqref{eq:flatflux}.

Consider the region between two compact globally regular levels, or
between such a level and a large enclosing sphere. By \textup{(P5)}
the latter boundary is noncritical. When the closure lies in $E$,
mollify the extension in \textup{(P4)} in a neighborhood of the
closure. Each mollification has zero divergence and converges uniformly
to the classical field on both boundaries, where \textup{(P2)} supplies
a noncritical collar. The divergence theorem therefore identifies
the two outward fluxes (with normals pointing away from $\Omega$).
The same equality holds at $\partial\Omega$ by taking regular levels
$u\uparrow1$ in the finite boundary collar. Such levels stay in the
collar by \textup{(P1)} and converge in $C^2$ by \textup{(P2)}.

Write $r=|x|$ and $e_r=x/r$. By \textup{(P5)}, uniformly in direction,
\[
 \begin{aligned}
 Dz&=\alpha_p\mathsf a r^{-\alpha_p-1}(e_r+o(1)),\\
 D^2z&=\alpha_p\mathsf a r^{-\alpha_p-2}
        [I-(\alpha_p+2)e_r\otimes e_r+o(1)].
 \end{aligned}
\]
The limiting tangential block is the identity. Homogeneity of
$T_{k-1}$ and $(\alpha_p+1)(p-1)=n-k$ give
\[
 r^{n-1}F_z\cdot e_r
 =\binom{n-1}{k-1}(\alpha_p\mathsf a)^{p-1}+o(1).
\]
The normalized gradient stays away from zero, so the error is uniform
also when the flux contains a negative gradient power. Integrating on
the sphere and dividing by
$\omega_{n-1}\binom{n-1}{k-1}\alpha_p^{p-1}$ in
\eqref{eq:capacitydef} yields \eqref{eq:Ccoefficient}.
\end{proof}

\subsection{The weighted Hessian}\label{sec:weightedhessian}

At a noncritical point the outward level normal and the Newton
flux with this orientation are
\[
 \nu=-\frac{Du}{|Du|},\qquad
 F_z=-|Du|^{p-k-1}T_{k-1}(-D^2u)Du.
\]
Since \(T_{k-1}(-D^2u)=(-1)^{k-1}T_{k-1}(D^2u)\),
\[
 \Div F_z=(-1)^k\LL u.
\]
Thus \(\mathcal L_{k,p}u=0\) is equivalent to \(\Div F_z=0\).
This choice of flux preserves the outward orientation used in
Assumption~\ref{ass:potential}\textup{(P4)}.

The divergence-free property of the Euclidean Newton tensor and
the identities
\[
 \tr\!\left(T_{k-1}(-D^2u)(-D^2u)\right)
   =k\sigma_k(-D^2u),\qquad
 D|Du|=-D^2u\,\nu.
\]
give
\[
 \begin{aligned}
 \Div F_z
 &=-|Du|^{p-k-1}\tr\!\left(T_{k-1}(-D^2u)D^2u\right)\\
 &\quad-(p-k-1)|Du|^{p-k-2}
       \ip{D|Du|}{T_{k-1}(-D^2u)Du}\\
 &=|Du|^{p-k-1}\Bigl[k\sigma_k(-D^2u)\\
 &\hspace{32mm}
       +(p-k-1)\ip{(-D^2u)T_{k-1}(-D^2u)\nu}{\nu}\Bigr].
 \end{aligned}
\]
Choose an orthonormal frame whose last vector is \(\nu\).
The sum of the principal \(k\)-minors involving that last index is
\[
 \sigma_k(-D^2u)-(T_k(-D^2u))_{\nu\nu}
 =\ip{(-D^2u)T_{k-1}(-D^2u)\nu}{\nu}.
\]
Indeed, the normal entry of \(T_k\) is the sum of the \(k\)-minors
avoiding that index, and the equality follows from the Newton
recurrence. In the expanded divergence the minors avoiding the
normal index therefore have coefficient \(k\), whereas those
involving it have coefficient \(k+(p-k-1)=p-1\).

Define
\[
 P=I+\left(\sqrt{\frac{p-1}{k}}-1\right)\nu\otimes\nu.
\]
It fixes tangent vectors and multiplies \(\nu\) by
\(\sqrt{(p-1)/k}\). Multiplying \(-D^2u\) by \(P\) on both sides
therefore multiplies each minor involving the normal index by
\((p-1)/k\) and leaves the other minors unchanged. Consequently,
\[
 \Div F_z=k|Du|^{p-k-1}\sigma_k(-PD^2uP).
\]
The weighted Hessian relevant to the admissible branch is thus
\(-PD^2uP\). Membership in \(\overline{\Gamma_k}\) does not
by itself imply strict ellipticity.

The same identity can be written without introducing another
parameter:
\[
 \Div F_z
 =k\sigma_k\!\left(D\!\left(-|Du|^{(p-k-1)/k}Du\right)\right).
\]
The Jacobian inside \(\sigma_k\) is
\(-|Du|^{(p-k-1)/k}P^2D^2u\).
Conjugation by \(P\) shows that \(-P^2D^2u\) is similar to
\(-PD^2uP\), and homogeneity supplies the factor
\(|Du|^{p-k-1}\). The Jacobian need not itself be symmetric.
When \(k=1\), the equation is the \(p\)-Laplace equation;
when \(p=k+1\), \(P=I\) and the weighted Hessian reduces to
\(-D^2u\).

\subsection{The conformal setting}\label{sec:conformal}

Set
\begin{equation}\label{eq:logvariables}
 \psi=-\frac{p-1}{n-k-p+1}\log u=-\frac{\log u}{\alpha_p}.
\end{equation}
Then $\psi=0$ on $\partial\Omega$, $\psi>0$ in $E$,
and $\psi\to\infty$ at infinity. Thus every finite closed slab of
$\psi$ in $\overline E$ is compact. Define
\begin{equation}\label{eq:conformal}
 \begin{gathered}
 g:=e^{-2\psi}\delta=u^{2(p-1)/(n-k-p+1)}\delta,\qquad
 \xi:=\nabla_g\psi,\\
 \KK:=g^{-1}\nabla_g^2\psi+|\nabla_g\psi|_g^2I.
 \end{gathered}
\end{equation}
All unmarked gradients, divergences, inner products, and tensor
contractions from now on use $g$. Euclidean operations retain a
subscript $\delta$ or the notation $D$.

The Hessian $\nabla_g^2\psi$ is the Hessian w.r.t $g$, whereas a
Newton transformation acts on an endomorphism. This explains the
factor $g^{-1}$ in \eqref{eq:conformal}: in coordinates,
$$
 (g^{-1}\nabla_g^2\psi)^i{}_j
 =g^{i\ell}(\nabla_g^2\psi)_{\ell j}.
$$
It raises the first index. In a $g$-orthonormal frame the numerical
matrix is the usual symmetric Hessian matrix. All eigenvalues of
$\KK$ and all tensors $T_j(\KK)$ refer to this self-adjoint
endomorphism. The metric is $C^1$ globally and $C^3$ on
$\mathcal R$; all local differential identities below are proved
on $\mathcal R$, where $|\nabla\psi|>0$.

\begin{lemma}
For any symmetric matrix $A$, vector $\zeta$, and scalar $t$,
\begin{equation}\label{eq:rankone}
 T_j(A+t\zeta\otimes\zeta)\zeta=T_j(A)\zeta.
\end{equation}
Here the outer product uses the inner product underlying the chosen
orthonormal frame.
\end{lemma}
\begin{proof}
 To prove the identity, use the generating polynomial
$$
 \det(I+sA)(I+sA)^{-1}=\sum_{j=0}^{n-1}s^jT_j(A).
$$
For sufficiently small $s$, put $M=I+sA$ and
$\beta=\zeta^{\mathsf T}M^{-1}\zeta$. The determinant lemma and
the rank-one inverse formula give
$$
 \det(M+st\zeta\zeta^{\mathsf T})
 =\det M\,(1+st\beta),\qquad
 (M+st\zeta\zeta^{\mathsf T})^{-1}\zeta
 =\frac{M^{-1}\zeta}{1+st\beta}.
$$
Multiplying cancels the two scalar factors. The resulting identity
between the two adjugate polynomials applied to $\zeta$ holds for
all small $s$, hence their coefficients agree. This proves
\eqref{eq:rankone} without any invertibility assumption on $A$.
\end{proof}

\begin{lemma}
On $\mathcal R$, where $|\nabla\psi|>0$, let
$$
 J_\delta=|D\psi|^{p-k-1}T_{k-1}(D^2\psi)D\psi.
$$
Then
\begin{equation}\label{eq:weightedflat}
 \Div_\delta(e^{-(n-k-p+1)\psi}J_\delta)=0.
\end{equation}
\end{lemma}
\begin{proof}
Differentiating $u=\exp\left(-\frac{n-k-p+1}{p-1}\psi\right)$ gives
$$
 \begin{aligned}
 Du&=-\frac{n-k-p+1}{p-1}uD\psi,\\
 D^2u&=-\frac{n-k-p+1}{p-1}u
       \left(D^2\psi-\frac{n-k-p+1}{p-1}D\psi\otimes D\psi\right).
 \end{aligned}
$$
The Newton tensor is homogeneous of degree $k-1$, including under
negative scalar multiplication. Consequently,
\begin{align*}
 &|Du|^{p-k-1}T_{k-1}(D^2u)Du\\
 &\quad=(-1)^k\left(\frac{n-k-p+1}{p-1}u\right)^{p-1}
       |D\psi|^{p-k-1}\\
 &\qquad\quad{}\times
 T_{k-1}\left(D^2\psi-\frac{n-k-p+1}{p-1}
                         D\psi\otimes D\psi\right)D\psi\\
 &\quad=(-1)^k\left(\frac{n-k-p+1}{p-1}u\right)^{p-1}J_\delta.
\end{align*}
The powers from the gradient weight, Newton tensor, and final gradient
add to $(p-k-1)+(k-1)+1=p-1$.
The second equality uses \eqref{eq:rankone}.
Since $u^{p-1}=e^{-(n-k-p+1)\psi}$, the last scalar factor is
$$
 (-1)^k\left(\frac{n-k-p+1}{p-1}\right)^{p-1}
 e^{-(n-k-p+1)\psi}.
$$
Thus the pointwise equation on $\mathcal R$ is equivalent to
\begin{equation}
 \Div_\delta(e^{-(n-k-p+1)\psi}J_\delta)=0.
\end{equation}
\end{proof}

Recall from \eqref{eq:conformal} that
$$
 \KK=g^{-1}\nabla_g^2\psi+|\nabla\psi|^2I,
$$
where $g^{-1}\nabla_g^2\psi$ is the Hessian with one index raised,
viewed as a self-adjoint endomorphism with respect to $g$.

\begin{lemma}
Let $J:=|\nabla\psi|^{p-k-1}T_{k-1}(\KK)\nabla_g\psi.$ Then on the noncritical set $\mathcal R$,
\begin{equation}\label{eq:conformalPDE}
 \Div_gJ=0.\qquad
\end{equation}
    \end{lemma}
\begin{proof}
    For $g=e^{-2\psi}\delta$, the connection formula reads
$$
 \nabla^g_VW=D_VW-V(\psi)W-W(\psi)V
                  +\delta(V,W)D\psi.
$$
Applying the definition of the Hessian therefore gives
$$
 \nabla_g^2\psi
 =D^2_\delta\psi+2\,\mathrm{d}\psi\otimes \mathrm{d}\psi
                     -|D\psi|_\delta^2\delta.
$$
Also $\xi=e^{2\psi}D\psi$ and $|\nabla\psi|=e^\psi|D\psi|_\delta$, so
$|\nabla\psi|^2g=|D\psi|_\delta^2\delta$. It follows that
$$
 D^2_\delta\psi
 =\nabla_g^2\psi-2\,\mathrm{d}\psi\otimes \mathrm{d}\psi+|\nabla\psi|^2g,
 \qquad
 \KK^i{}_j=e^{2\psi}
   \bigl(\psi_{ij}+2\psi_i\psi_j\bigr)
$$
in Euclidean coordinates, with the latter displayed array understood
as a mixed-index matrix. Applying homogeneity and
\eqref{eq:rankone} once more gives
\begin{align*}
 J&:=|\nabla\psi|^{p-k-1}T_{k-1}(\KK)\xi\\
 &=e^{(p-k-1)\psi}e^{2(k-1)\psi}e^{2\psi}
      |D\psi|^{p-k-1}
      T_{k-1}(D^2\psi+2D\psi\otimes D\psi)D\psi\\
 &=e^{(p+k-1)\psi}J_\delta.
\end{align*}
The exponent is $(p-k-1)+2(k-1)+2=p+k-1$. On the other hand,
$\mathrm{d}V_g=e^{-n\psi}\mathrm{d}V_\delta$ and
$n-(p+k-1)=n-k-(p-1)=(n-k-p+1)$. The coordinate formula for divergence
therefore yields
$$
 \Div_gJ=e^{n\psi}\Div_\delta(e^{-n\psi}J)
        =e^{n\psi}\Div_\delta(e^{-(n-k-p+1)\psi}J_\delta).
$$
The equation on the noncritical set is therefore \eqref{eq:conformalPDE}.
For $k=1$ this is $\Delta_p^g\psi=0$. For $k>1$ the equation
retains the shifted tensor $\KK$; the term $|\nabla\psi|^2I$ is part of the
exact transformation.
\end{proof}

\begin{lemma}
    Let $\HH:=h+|\nabla\psi|I$ on that tangent space. Then on $\mathcal R$
\begin{equation}\label{eq:leveltransform}
 h_\delta=e^{-\psi}\HH,\qquad
 |D\psi|_\delta=e^{-\psi}|\nabla\psi|,\qquad
 \mathrm{d}S_\delta=e^{(n-1)\psi}\mathrm{d}S_g.
\end{equation}
\end{lemma}
\begin{proof}
    On a local level set let $\nu=\xi/|\nabla\psi|$ and let $h$ be the
$g$-shape operator. Its tangent space has dimension $n-1$.
 The Euclidean normal with
the same outward orientation is
$\nu_\delta=D\psi/|D\psi|=-Du/|Du|$, and
$\nu=e^\psi\nu_\delta$. If $V$ is tangent to a level, then
$V(\psi)=0$ and $\delta(V,\nu_\delta)=0$. Thus
$$
 \nabla^g_V\nu=e^\psi D_V\nu_\delta
       -e^\psi\nu_\delta(\psi)V
       =e^\psi h_\delta(V)-|\nabla\psi|V.
$$
This proves $h=e^\psi h_\delta-|\nabla\psi|I$ and hence
\begin{equation}
 h_\delta=e^{-\psi}\HH,\qquad
 |D\psi|_\delta=e^{-\psi}|\nabla\psi|,\qquad
 \mathrm{d}S_\delta=e^{(n-1)\psi}\mathrm{d}S_g.
\end{equation}
In particular,
$\sigma_{k-1}(\HH)=e^{(k-1)\psi}\sigma_{k-1}(h_\delta)>0$
by \textup{(P3)}.
\end{proof}

\begin{proposition}[Orthogonal decomposition and the higher-order Kato replacement]
\label{prop:kato}
On $\mathcal R$, recall that
$$
 \begin{gathered}
 \nu=\frac{\nabla_g\psi}{|\nabla\psi|},\qquad
 \HH=h+|\nabla\psi|I=e^\psi h_\delta,\\
 \KK=g^{-1}\nabla_g^2\psi+|\nabla\psi|^2I.
 \end{gathered}
$$
Here $\nabla^\top |\nabla\psi|$ is the gradient of $|\nabla\psi|$ tangent to the level
set, and $\nu(|\nabla\psi|)$ is its derivative along the unit normal.
All inner products, norms, and tensor contractions below use $g$.
By \textup{(P3)}, $\sigma_{k-1}(\HH)>0$.

\emph{(i) Tangential and normal decomposition.}
The gradient and Hessian split as
\begin{equation}\label{eq:orthogonal}
 \begin{aligned}
 \nabla |\nabla\psi|&=\nabla^\top |\nabla\psi|+\nu(|\nabla\psi|)\nu,\\
 \bigl|\nabla|\nabla\psi|\bigr|^2&=\bigl|\nabla^\top|\nabla\psi|\bigr|^2+\nu(|\nabla\psi|)^2,\\
 |\nabla_g^2\psi|^2
 &=|\nabla\psi|^2|h|^2+2\bigl|\nabla^\top|\nabla\psi|\bigr|^2+\nu(|\nabla\psi|)^2.
 \end{aligned}
\end{equation}
In a $g$-orthonormal frame adapted to the level set,
\begin{equation}\label{eq:blockK}
 \KK=\begin{pmatrix}
 |\nabla\psi|\HH&\nabla^\top |\nabla\psi|\\
 (\nabla^\top |\nabla\psi|)^{\mathsf T}&\nu(|\nabla\psi|)+|\nabla\psi|^2
 \end{pmatrix}.
\end{equation}
The tangent vector in this matrix is represented by its $n-1$
components in the chosen frame.

\emph{(ii) Newton tensors and their divergence.}
For $0\le j\le k$, expansion of the principal minors gives
\begin{equation}\label{eq:blockminor}
 \begin{aligned}
 \sigma_{j+1}(\KK)
 &=|\nabla\psi|^{j+1}\sigma_{j+1}(\HH)
   +(\nu(|\nabla\psi|)+|\nabla\psi|^2)|\nabla\psi|^j\sigma_j(\HH)\\
 &\quad-|\nabla\psi|^{j-1}T_{j-1}(\HH)
                  (\nabla^\top |\nabla\psi|,\nabla^\top |\nabla\psi|).
 \end{aligned}
\end{equation}
In particular, the normal and mixed components are
\begin{equation}\label{eq:blocknewton}
 \begin{aligned}
 (T_j(\KK))_{\nu\nu}&=|\nabla\psi|^j\sigma_j(\HH),\\
 (T_j(\KK))_{\top\nu}
 &=-|\nabla\psi|^{j-1}T_{j-1}(\HH)\nabla^\top |\nabla\psi|.
 \end{aligned}
\end{equation}
Terms involving $T_{-1}$ are zero. The divergence in the conformal
metric is
\begin{equation}\label{eq:newtondiv}
 \Div_gT_j(\KK)
 =\bigl[(n-j)\sigma_j(\KK)I-nT_j(\KK)\bigr]\nabla_g\psi.
\end{equation}

\emph{(iii) The equation in level-set quantities.}
The local equation \eqref{eq:conformalPDE} is equivalent to
\begin{equation}\label{eq:levelPDE}
 \begin{aligned}
 &(p-1)\left[\nu(|\nabla\psi|)
  -\frac{T_{k-2}(\HH)(\nabla^\top |\nabla\psi|,\nabla^\top |\nabla\psi|)}{|\nabla\psi|\,\sigma_{k-1}(\HH)}\right]\\
 &\qquad=|\nabla\psi|\left[(n-k)|\nabla\psi|
       -k\frac{\sigma_k(\HH)}{\sigma_{k-1}(\HH)}\right].
 \end{aligned}
\end{equation}
The quadratic term on the left is the contribution of the mixed
Hessian block. It vanishes when $k=1$ because $T_{-1}=0$.

\emph{(iv) Nonnegative curvature expressions.}
The following inequalities hold, with the first understood as
an inequality of quadratic forms on the tangent space:
\begin{equation}\label{eq:BE}
 \begin{aligned}
 T_{k-1}(\HH)
  -\frac{\sigma_k(\HH)}{\sigma_{k-1}(\HH)}T_{k-2}(\HH)&\ge0,\\
 \frac{k(n-k-1)}{n-k}
  \frac{\sigma_k(\HH)^2}{\sigma_{k-1}(\HH)}
  -(k+1)\sigma_{k+1}(\HH)&\ge0.
 \end{aligned}
\end{equation}
They are related by the weighted curvature identity
\begin{equation}\label{eq:weightedcurvature}
 \begin{aligned}
 &\tr\!\left[
  \left(T_{k-1}(\HH)
   -\frac{\sigma_k(\HH)}{\sigma_{k-1}(\HH)}T_{k-2}(\HH)\right)
   \HH^2\right]
   -\frac{k}{n-k}\frac{\sigma_k(\HH)^2}{\sigma_{k-1}(\HH)}\\
 &\qquad=\frac{k(n-k-1)}{n-k}
   \frac{\sigma_k(\HH)^2}{\sigma_{k-1}(\HH)}
   -(k+1)\sigma_{k+1}(\HH).
 \end{aligned}
\end{equation}
These formulas express the curvature terms directly through the
elementary symmetric polynomials and Newton tensors of $\HH$.
\end{proposition}

\begin{proof}
\emph{Level-set blocks.}
Choose a $g$-orthonormal frame $e_1,\ldots,e_{n-1},\nu$.
For tangent vectors $V,W$, differentiating $\nabla_g\psi=|\nabla\psi|\nu$ gives
$$
 \nabla_g^2\psi(V,W)
 =g(\nabla_V(|\nabla\psi|\nu),W)=|\nabla\psi|\,g(hV,W).
$$
The mixed and normal entries are
$$
 \nabla_g^2\psi(V,\nu)=V(|\nabla\psi|),\qquad
 \nabla_g^2\psi(\nu,\nu)=\nu(|\nabla\psi|).
$$
Thus the Hessian matrix in this frame is
$$
 \begin{pmatrix}
 |\nabla\psi|h&\nabla^\top |\nabla\psi|\\
 (\nabla^\top |\nabla\psi|)^{\mathsf T}&\nu(|\nabla\psi|)
 \end{pmatrix}.
$$
The squared tensor norm sums all matrix entries, including both
copies of the mixed block. This proves \eqref{eq:orthogonal}.
Adding $|\nabla\psi|^2I$ to the entire $n$-dimensional endomorphism gives
\eqref{eq:blockK} because $|\nabla\psi|h+|\nabla\psi|^2I=|\nabla\psi|\HH$ on the tangent space.

\emph{Newton minors and components.}
Diagonalize $\HH$ with eigenvalues $\kappa_1,\ldots,\kappa_{n-1}$
in the tangent frame. A principal $(j+1)$-minor not containing
the normal direction is a product of $j+1$ entries of $|\nabla\psi|\HH$.
Their sum is $|\nabla\psi|^{j+1}\sigma_{j+1}(\HH)$.
A minor containing the normal direction and tangent indices
$I$, with $|I|=j$, is
$$
 (\nu(|\nabla\psi|)+|\nabla\psi|^2)|\nabla\psi|^j\prod_{i\in I}\kappa_i
 -|\nabla\psi|^{j-1}\sum_{i\in I}e_i(|\nabla\psi|)^2
          \prod_{\ell\in I\setminus\{i\}}\kappa_\ell.
$$
This determinant expansion divides by no eigenvalue, so it also
holds when some $\kappa_i$ vanish. Write $\kappa|i$ for the list
with $\kappa_i$ omitted. For fixed $i$, summing the products in
the second term over the remaining $j-1$ indices gives
$\sigma_{j-1}(\kappa|i)$, the $i$th eigenvalue of
$T_{j-1}(\HH)$. Since $e_i(|\nabla\psi|)$ is the $i$th component of
$\nabla^\top |\nabla\psi|$, summing all minors gives \eqref{eq:blockminor}.

The differential identity
$\mathrm{d}\sigma_{j+1}(\KK)[\dot{\KK}]
=\tr(T_j(\KK)\dot{\KK})$
shows that differentiating the minor formula with respect to
$\KK_{\nu\nu}=\nu(|\nabla\psi|)+|\nabla\psi|^2$ gives
$(T_j(\KK))_{\nu\nu}=|\nabla\psi|^j\sigma_j(\HH)$.
Varying the mixed entry $e_i(|\nabla\psi|)$ changes both $\KK_{i\nu}$ and
$\KK_{\nu i}$, so its derivative is $2(T_j(\KK))_{i\nu}$.
The quadratic term differentiates to
$-2|\nabla\psi|^{j-1}(T_{j-1}(\HH)\nabla^\top |\nabla\psi|)_i$.
Dividing by two proves \eqref{eq:blocknewton}. In particular,
$$
 \begin{aligned}
 \langle J,\nu\rangle_g
 &=|\nabla\psi|^{p-k-1}|\nabla\psi|\,(T_{k-1}(\KK))_{\nu\nu}\\
 &=|\nabla\psi|^{p-k-1+1+k-1}\sigma_{k-1}(\HH)
  =|\nabla\psi|^{p-1}\sigma_{k-1}(\HH).
 \end{aligned}
$$
Every inner product in this calculation uses $g$.

\emph{Divergence of the Newton tensor.}
Set $\rho=e^{2\psi}/2$. Direct differentiation gives
$$
 (D^2_\delta\rho)_{ij}
 =e^{2\psi}(\psi_{ij}+2\psi_i\psi_j)=\KK^i{}_j
$$
as numerical arrays in Euclidean coordinates. Thus $\KK$ is the
Euclidean Hessian endomorphism of $\rho$ in these coordinates.
For $j\ge1$, the entries of $T_j(\KK)$ have the generalized
Kronecker formula
$$
 (T_j)^i{}_r
 =\frac{1}{j!}
 \delta^{i i_1\cdots i_j}_{r r_1\cdots r_j}
 \rho_{i_1r_1}\cdots\rho_{i_jr_j}.
$$
After applying $\partial_i$, a derivative landing on the $a$th
factor contains $\rho_{i i_a r_a}$, which is symmetric in
$i,i_a$. The Kronecker tensor is antisymmetric in the same two
upper indices, so that summand vanishes after contraction.
Therefore $\partial_i(T_j)^i{}_r=0$; the same statement is immediate
for $T_0=I$.

This Euclidean divergence identity does not itself say that the
$g$-divergence is zero. The conformal Christoffel symbols are
$$
 \Gamma^\ell_{ir}
 =-\delta^\ell_i\psi_r-\delta^\ell_r\psi_i
                         +\delta_{ir}\psi_\ell,
 \qquad \Gamma^i_{i\ell}=-n\psi_\ell.
$$
The matrix $T_j(\KK)$ is symmetric also as a Euclidean numerical
matrix. The divergence of this mixed tensor is therefore
\begin{align*}
 \nabla_i(T_j)^i{}_r
 &=\partial_i(T_j)^i{}_r+\Gamma^i_{i\ell}(T_j)^\ell{}_r
                               -\Gamma^\ell_{ir}(T_j)^i{}_\ell\\
 &=-n\psi_\ell(T_j)^\ell{}_r
   +(\tr T_j)\psi_r+\psi_i(T_j)^i{}_r-\psi_\ell(T_j)^r{}_\ell\\
 &=(\tr T_j)\psi_r-n(T_j)^i{}_r\psi_i.
\end{align*}
The last two terms on the second line cancel by symmetry.
Raising the free index $r$ with $g^{rs}$ turns $\mathrm{d}\psi$ into
$\nabla_g\psi$. The trace identity
$\tr T_j(\KK)=(n-j)\sigma_j(\KK)$ proves \eqref{eq:newtondiv}.

\emph{The level equation.}
Recall that $\xi=\nabla_g\psi=|\nabla\psi|\nu$ and
$\nabla\xi=\KK-|\nabla\psi|^2I$. The product rule gives
\begin{align*}
 \Div_g\bigl(T_j(\KK)\xi\bigr)
 &=\langle\Div_gT_j(\KK),\xi\rangle_g
                      +\tr\bigl(T_j(\KK)(\KK-|\nabla\psi|^2I)\bigr)\\
 &=(n-j)\sigma_j(\KK)|\nabla\psi|^2
   -n\langle T_j(\KK)\xi,\xi\rangle_g\\
 &\quad +(j+1)\sigma_{j+1}(\KK)-|\nabla\psi|^2\tr T_j(\KK)\\
 &=(j+1)\sigma_{j+1}(\KK)
   -n\langle T_j(\KK)\xi,\xi\rangle_g.
\end{align*}
The cancellation uses
$\tr(T_j(\KK)\KK)=(j+1)\sigma_{j+1}(\KK)$ and
$\tr T_j(\KK)=(n-j)\sigma_j(\KK)$.
For any constants $\mu,\alpha$, differentiating the scalar weight
$e^{-\mu\psi}|\nabla\psi|^\alpha$ gives the useful general formula
\begin{equation}\label{eq:universaldiv}
 \begin{aligned}
 &\Div_g\bigl(e^{-\mu\psi}|\nabla\psi|^\alpha T_j(\KK)\xi\bigr)\\
 &\quad=e^{-\mu\psi}|\nabla\psi|^\alpha\Bigl[
    (j+1)\sigma_{j+1}(\KK)
    -(n+\mu)\langle T_j(\KK)\xi,\xi\rangle_g\\
 &\hspace{42mm}
    +\frac{\alpha}{|\nabla\psi|}\langle T_j(\KK)\xi,\nabla |\nabla\psi|\rangle_g\Bigr].
 \end{aligned}
\end{equation}
Both contractions are explicit consequences of \eqref{eq:blocknewton}:
$$
 \begin{aligned}
 \langle T_j(\KK)\xi,\xi\rangle_g
 &=|\nabla\psi|^{j+2}\sigma_j(\HH),\\
 \langle T_j(\KK)\xi,\nabla |\nabla\psi|\rangle_g
 &=|\nabla\psi|^{j+1}\sigma_j(\HH)\nu(|\nabla\psi|)
   -|\nabla\psi|^jT_{j-1}(\HH)(\nabla^\top |\nabla\psi|,\nabla^\top |\nabla\psi|).
 \end{aligned}
$$
Indeed, multiply both components in \eqref{eq:blocknewton} by
$|\nabla\psi|$, since $\xi=|\nabla\psi|\nu$, and pair them with
the decomposition of $\nabla|\nabla\psi|$ in \eqref{eq:orthogonal}.

Apply \eqref{eq:universaldiv} with $\mu=0$,
$\alpha=p-k-1$, and $j=k-1$. Before collecting terms, this is
\begin{align*}
 |\nabla\psi|^{-(p-k-1)}\Div_gJ
 &=k|\nabla\psi|^k\sigma_k(\HH)
    +k(\nu(|\nabla\psi|)+|\nabla\psi|^2)|\nabla\psi|^{k-1}\sigma_{k-1}(\HH)\\
 &\quad-k|\nabla\psi|^{k-2}T_{k-2}(\HH)(\nabla^\top |\nabla\psi|,\nabla^\top |\nabla\psi|)
    -n|\nabla\psi|^{k+1}\sigma_{k-1}(\HH)\\
 &\quad+(p-k-1)|\nabla\psi|^{k-1}\sigma_{k-1}(\HH)\nu(|\nabla\psi|)\\
 &\quad-(p-k-1)|\nabla\psi|^{k-2}
            T_{k-2}(\HH)(\nabla^\top |\nabla\psi|,\nabla^\top |\nabla\psi|).
\end{align*}
The two terms containing $\nu(|\nabla\psi|)$ combine with coefficient $p-1$,
as do the two mixed-block contractions. Thus
$$
 \begin{aligned}
 \Div_gJ=|\nabla\psi|^{p-3}\Bigl\{
 &(p-1)|\nabla\psi|\,\sigma_{k-1}(\HH)\nu(|\nabla\psi|)\\
 &-(p-1)T_{k-2}(\HH)(\nabla^\top |\nabla\psi|,\nabla^\top |\nabla\psi|)\\
 &+k|\nabla\psi|^2\sigma_k(\HH)-(n-k)|\nabla\psi|^3\sigma_{k-1}(\HH)\Bigr\}.
 \end{aligned}
$$
Since $|\nabla\psi|>0$ and $\sigma_{k-1}(\HH)>0$, setting this divergence
equal to zero and dividing its bracket by
$|\nabla\psi|\,\sigma_{k-1}(\HH)$ proves \eqref{eq:levelPDE}.

To prove
\eqref{eq:BE} and \eqref{eq:weightedcurvature}, we use Newton's inequality for real eigenvalues.
The curvature inequalities needed here do not require every
principal curvature to be positive. For any list
$\lambda=(\lambda_1,\ldots,\lambda_m)$ of real numbers,
\begin{equation}\label{eq:realNewton}
 \sigma_j(\lambda)^2\ge
 \frac{(j+1)(m-j+1)}{j(m-j)}
 \sigma_{j-1}(\lambda)\sigma_{j+1}(\lambda),
 \qquad 1\le j\le m-1.
\end{equation}
To verify the coefficient, set
$$
 P(t)=\prod_{i=1}^m(t+\lambda_i)
     =\sum_{\ell=0}^m\sigma_\ell(\lambda)t^{m-\ell},
 \qquad f=P^{(m-j-1)}.
$$
Repeated Rolle's theorem shows that $f$ has $j+1$ real roots,
counted with multiplicity; denote them by $r_1,\ldots,r_{j+1}$.
Away from these roots, logarithmic differentiation gives
$$
 \frac{f'}{f}=\sum_{i=1}^{j+1}\frac{1}{t-r_i},\qquad
 \frac{f''}{f}=
  \left(\sum_{i=1}^{j+1}\frac{1}{t-r_i}\right)^2
     -\sum_{i=1}^{j+1}\frac{1}{(t-r_i)^2}.
$$
Consequently,
$$
 (f')^2-\frac{j+1}{j}ff''
 =\frac{f^2}{j}\left[
  (j+1)\sum_{i=1}^{j+1}\frac{1}{(t-r_i)^2}
   -\left(\sum_{i=1}^{j+1}\frac{1}{t-r_i}\right)^2\right]\ge0
$$
by Cauchy--Schwarz. Polynomial continuity extends the inequality
to the roots of $f$. At $t=0$,
$$
 \begin{aligned}
 f(0)&=(m-j-1)!\,\sigma_{j+1}(\lambda),\\
 f'(0)&=(m-j)!\,\sigma_j(\lambda),\\
 f''(0)&=(m-j+1)!\,\sigma_{j-1}(\lambda).
 \end{aligned}
$$
Substitution and division by $((m-j)!)^2$ give the factor
$$
 \frac{j+1}{j}
 \frac{(m-j-1)!(m-j+1)!}{((m-j)!)^2}
 =\frac{(j+1)(m-j+1)}{j(m-j)},
$$
which proves \eqref{eq:realNewton}. We also record the equality
criterion used below. If $\sigma_{j+1}>0$, then $f(0)\ne0$ and equality
in Cauchy--Schwarz at $t=0$ makes all $r_i$ equal. Comparing the
first two normalized coefficients of the monic polynomials proportional
to $P$ and to $f$ gives
$(\sigma_1/m)^2=\sigma_2/\binom m2$ for the original eigenvalues.
Equivalently, $\sum_{a<b}(\lambda_a-\lambda_b)^2=0$. Thus they all coincide.

\emph{The signs of the two curvature terms.}
Apply \eqref{eq:realNewton} to the $n-1$ real eigenvalues of
$\HH$, with index $j=k$. The hypothesis $1\le k\le n-2$
puts this index in the required range and gives
$$
 \sigma_k(\HH)^2\ge
 \frac{(k+1)(n-k)}{k(n-k-1)}
 \sigma_{k-1}(\HH)\sigma_{k+1}(\HH).
$$
Since $\sigma_{k-1}(\HH)>0$, multiplication by
$k(n-k-1)/(n-k)$ and division by $\sigma_{k-1}(\HH)$ yield
$$
 \frac{k(n-k-1)}{n-k}
 \frac{\sigma_k(\HH)^2}{\sigma_{k-1}(\HH)}
 -(k+1)\sigma_{k+1}(\HH)\ge0.
$$

For the tensor inequality, choose an orthonormal eigenbasis
$e_1,\ldots,e_{n-1}$ of $\HH$, with eigenvalues
$\kappa=(\kappa_1,\ldots,\kappa_{n-1})$.
Both Newton tensors in \eqref{eq:BE} are diagonal in this basis.
For $k\ge2$, write $\kappa|i$ for the $n-2$ eigenvalues
remaining after $\kappa_i$ is omitted. The Newton eigenvalue
formula is
$$
 T_{k-1}(\HH)e_i=\sigma_{k-1}(\kappa|i)e_i,
 \qquad
 T_{k-2}(\HH)e_i=\sigma_{k-2}(\kappa|i)e_i.
$$
Separating products that contain $\kappa_i$ from those that do not
gives
\begin{align*}
 \sigma_{k-1}(\HH)
 &=\sigma_{k-1}(\kappa|i)
      +\kappa_i\sigma_{k-2}(\kappa|i),\\
 \sigma_k(\HH)
 &=\sigma_k(\kappa|i)
      +\kappa_i\sigma_{k-1}(\kappa|i).
\end{align*}
Consequently, the $i$th eigenvalue of the tensor in
\eqref{eq:BE}, multiplied by $\sigma_{k-1}(\HH)$, is
\begin{align*}
 &\sigma_{k-1}(\HH)\sigma_{k-1}(\kappa|i)
       -\sigma_k(\HH)\sigma_{k-2}(\kappa|i)\\
 &\qquad=\sigma_{k-1}(\kappa|i)^2
       -\sigma_k(\kappa|i)\sigma_{k-2}(\kappa|i).
\end{align*}
If the product in the last term is nonpositive, this expression
is nonnegative directly. If that product is positive, apply
\eqref{eq:realNewton} to the $n-2$ numbers $\kappa|i$, now with
$j=k-1$. The range $2\le k\le n-2$ ensures
$1\le j\le(n-2)-1$, so
$$
 \sigma_{k-1}(\kappa|i)^2
 \ge\frac{k(n-k)}{(k-1)(n-k-1)}
      \sigma_{k-2}(\kappa|i)\sigma_k(\kappa|i)
 \ge\sigma_{k-2}(\kappa|i)\sigma_k(\kappa|i).
$$
The second inequality uses both the positive product and the fact
that its coefficient is greater than one. Division by
$\sigma_{k-1}(\HH)>0$ proves that every eigenvalue of the tensor
in \eqref{eq:BE} is nonnegative. For $k=1$ that tensor is simply
$T_0-(\sigma_1/\sigma_0)T_{-1}=I$.
Thus neither curvature inequality requires every eigenvalue of
$\HH$ to be positive.
They both hold for any self-adjoint $\mathcal H$ with $\sigma_{k-1}(\mathcal H)>0$; no additional Gårding-cone assumption is needed.

\emph{The weighted trace identity.}
The Newton recurrence gives
$\HH T_{k-1}(\HH)=\sigma_k(\HH)I-T_k(\HH)$.
Multiply by $\HH$, take the trace, and use
$\tr(\HH T_k(\HH))=(k+1)\sigma_{k+1}(\HH)$ to obtain
$$
 \tr\bigl(T_{k-1}(\HH)\HH^2\bigr)
 =\sigma_1(\HH)\sigma_k(\HH)-(k+1)\sigma_{k+1}(\HH).
$$
The same calculation one index lower gives
$$
 \tr\bigl(T_{k-2}(\HH)\HH^2\bigr)
 =\sigma_1(\HH)\sigma_{k-1}(\HH)-k\sigma_k(\HH).
$$
For $k=1$ the second identity still holds: both sides are zero
because $T_{-1}=0$ and $\sigma_0=1$.

Subtract $\sigma_k(\HH)/\sigma_{k-1}(\HH)$ times the second
identity from the first. The terms containing $\sigma_1(\HH)$
cancel, leaving
\begin{align*}
 &\tr\!\left[
  \left(T_{k-1}(\HH)
      -\frac{\sigma_k(\HH)}{\sigma_{k-1}(\HH)}T_{k-2}(\HH)
  \right)\HH^2\right]\\
 &\qquad=
  k\frac{\sigma_k(\HH)^2}{\sigma_{k-1}(\HH)}
      -(k+1)\sigma_{k+1}(\HH).
\end{align*}
Subtracting
$k\sigma_k(\HH)^2/((n-k)\sigma_{k-1}(\HH))$
from both sides changes the coefficient $k$ on the right to
$k(n-k-1)/(n-k)$. This is precisely
\eqref{eq:weightedcurvature}, whose right-hand side is
nonnegative by \eqref{eq:BE}. \end{proof}

\begin{remark}
When $k=1$, $\sigma_0(\HH)=1$ and $T_{-1}(\HH)=0$.
Since $\sigma_1(\HH)=\tr h+(n-1)|\nabla\psi|$, the equation
\eqref{eq:levelPDE} reduces to
$$
 (p-1)\nu(|\nabla\psi|)+|\nabla\psi|\,\tr h=0.
$$
Write $h^\circ=h-(\tr h)I/(n-1)$ for the trace-free shape
operator. The usual trace decomposition and the preceding equation give
$$
 |h|^2=|h^\circ|^2+\frac{(\tr h)^2}{n-1},\qquad
 |\nabla\psi|^2|h|^2=|\nabla\psi|^2|h^\circ|^2+\frac{(p-1)^2}{n-1}\nu(|\nabla\psi|)^2.
$$
Inserting this into \eqref{eq:orthogonal} yields
$$
 |\nabla_g^2\psi|^2
 =|\nabla\psi|^2|h^\circ|^2+2\bigl|\nabla^\top|\nabla\psi|\bigr|^2
   +\left(1+\frac{(p-1)^2}{n-1}\right)\nu(|\nabla\psi|)^2.
$$
Subtracting
$\bigl(1+(p-1)^2/(n-1)\bigr)\bigl|\nabla|\nabla\psi|\bigr|^2$
and using the orthogonal decomposition of $\nabla |\nabla\psi|$ proves
\begin{equation}\label{eq:scalarKato}
 \begin{aligned}
 |\nabla_g^2\psi|^2
 -\left(1+\frac{(p-1)^2}{n-1}\right)\bigl|\nabla|\nabla\psi|\bigr|^2
 =|\nabla\psi|^2|h^\circ|^2
 +\left(1-\frac{(p-1)^2}{n-1}\right)\bigl|\nabla^\top|\nabla\psi|\bigr|^2.
 \end{aligned}
\end{equation}
This is an analogous identity of Kato's type (where $C=1+\frac{(p-1)^2}{n-1}>1$); the coefficient of $\bigl|\nabla^\top|\nabla\psi|\bigr|^2$
need not be positive throughout the parameter range.
The first expression in \eqref{eq:BE} is $I$ when $k=1$.
For the second, the identity
$2\sigma_2(\HH)=(\tr\HH)^2-|\HH|^2$ gives
$$
 \begin{aligned}
 \frac{n-2}{n-1}\sigma_1(\HH)^2-2\sigma_2(\HH)
 &=|\HH|^2-\frac{(\tr\HH)^2}{n-1}\\
 &=|h^\circ|^2.
 \end{aligned}
$$
For general $k$, Proposition~\ref{prop:kato} supplies the block
identities and nonnegative curvature expressions used in place
of this scalar identity. It is not an assertion of a universal
higher-order Kato theorem on an arbitrary Riemannian manifold.
\end{remark}

\begin{definition}[Conformal level functional]\label{def:Phi}
At a globally regular value $s\ge0$, set $\Sigma_s=\{\psi=s\}$,
with $\Sigma_0=\partial\Omega$, and define
\begin{equation}\label{eq:Phi}
 \Phi(s)=\int_{\Sigma_s}|\nabla\psi|^p\sigma_{k-1}(\HH)\dd S_g,
 \qquad \widehat\Phi(s)= \frac{\Phi(s)}{\omega_{n-1}\binom{n-1}{k-1}}.
\end{equation}
No value or differentiability of this geometric integral at a
critical level is imposed by this definition.
\end{definition}

\paragraph{The Euclidean level functional.}
On $\Sigma_s$, the conformal factors in
$|\nabla\psi|^p$, $\sigma_{k-1}(\HH)$, and $\dd S_g$ combine to give
\begin{equation}\label{eq:Phiflat}
 \Phi(s)=e^{-(n-k-p)s}
 \int_{\{\psi=s\}}|D\psi|^p\sigma_{k-1}(h_\delta)\dd S_\delta.
\end{equation}
In terms of the original potential the same functional is given by
\begin{equation}\label{eq:U}
 U_{k,p}(\tau)=\tau^{-\frac{n-k}{n-k-p+1}}
 \int_{\{u=\tau\}}|Du|^p\sigma_{k-1}(h_\delta)\dd S_\delta.
\end{equation}
Indeed, substituting $u=e^{-\alpha_p\psi}$ in \eqref{eq:Phiflat} gives
\begin{equation}\label{eq:Urelation}
 \Phi(s)=\alpha_p^{-p} U_{k,p}(e^{-\alpha_p s}).
\end{equation}
For $k=1$, the normalization $\varphi=(n-2)\psi$ used in
\cite{AFM} gives $\Phi_p^{\mathrm{AFM}}((n-2)s)=(n-2)^p\Phi(s)$.

\paragraph{Radial potentials.}
The conserved flux of a decreasing radial solution $f$ on the exterior
of $B_R$ is proportional to $r^{n-k}(-f')^{p-1}$ by
\eqref{eq:flatflux}. Integrating this first-order equation with boundary
values $1$ and $0$ yields
\begin{equation}\label{eq:ball}
 u=(R/r)^{\alpha_p},\qquad \psi=\log(r/R),\qquad
 g=R^2\bigl((\mathrm{d}\psi)^2+g_{\Sphere^{n-1}}\bigr).
\end{equation}
The slices have $g$-area $\omega_{n-1}R^{n-1}$ and $g$-shape operator
zero. Thus $|\nabla\psi|=R^{-1}$ and $\HH=R^{-1}I$, giving
$\widehat\Phi(s)=R^{n-k-p}$ and $\CC(B_R)=R^{n-k-p+1}$.

\paragraph{Asymptotic level geometry.}
For a general potential satisfying \textup{(P5)}, put
$R_\infty=\mathsf a^{\frac{p-1}{n-k-p+1}}$. Since
$\psi=-\frac{p-1}{n-k-p+1}\log u$, the first two derivatives are
$$
 D\psi=-\frac{p-1}{n-k-p+1}\frac{Du}{u},\qquad
 D^2\psi=-\frac{p-1}{n-k-p+1}
       \left(\frac{D^2u}{u}-\frac{Du\otimes Du}{u^2}\right).
$$
Using all three estimates in \textup{(P5)}, including the second-derivative
estimate, yields uniformly in angle
$$
 \psi=\log(r/R_\infty)+o(1),\quad
 D\psi=r^{-1}e_r+o(r^{-1}),\quad
 D^2\psi=r^{-2}(I-2e_r\otimes e_r)+o(r^{-2}).
$$
In particular $\partial_r\psi=r^{-1}(1+o(1))>0$ for large $r$.
Fix a radius beyond which this is true. Since $\psi$ is bounded
on the remaining compact region, every sufficiently large level
lies entirely beyond that radius. On each radial ray, $\psi$ then
increases strictly to infinity. It follows that the entire level
$\Sigma_s$ is a unique radial graph $r=r_s(\theta)$.

In logarithmic radial coordinates,
$F(t,\theta)=\psi(R_\infty e^t\theta)=t+o(1)$, with first and second
derivatives of the error also $o(1)$ in the cylinder metric. Indeed,
$F_t=r\partial_r\psi$, $F_{tt}=r\partial_r\psi+r^2\partial_{rr}\psi$,
and the spherical Hessian is
$r^2 D^2\psi|_{T\Sphere}-r\partial_r\psi\,g_{\Sphere}$; the preceding
expansions give the claim, including the mixed derivatives.
Since $F_t\to1$, twice differentiating
$F(\log(r_s/R_\infty),\theta)=s$ gives
\[
 f_s(\theta)=\frac{r_s(\theta)}{R_\infty e^s}
                 \longrightarrow1\quad\text{in }C^2(\Sphere^{n-1}).
\]
The induced conformal metric is
$R_\infty^2(f_s^2 g_{\Sphere}+\mathrm{d}f_s\otimes \mathrm{d}f_s)$, so its area density
tends uniformly to $R_\infty^{n-1}\dd S_{\Sphere}$. The corresponding
gradient length tends to $R_\infty^{-1}$ by the gradient expansion.
By $C^2$ convergence the Euclidean shape operator of the graph
$f_s(\theta)\theta$ tends to $I$; restoring its radius gives
$\HH=e^s h_\delta\to R_\infty^{-1}I$. Thus \eqref{eq:Phi} yields
\begin{equation}\label{eq:Phiinfinity}
 \widehat\Phi(\infty)
 =R_\infty^{n-k-p}
 =\mathsf a^{(n-k-p)/\alpha_p}
 =\CC(\Omega)^{(n-k-p)/(n-k-p+1)}.
\end{equation}
The geometric limit uses \textup{(P5)}. The last equality uses
Lemma~\ref{lem:flux}, which also requires \textup{(P4)}.

\section{The monotonicity}
\label{sec:proof}

In this section we prove Theorem~A, 
following the effective monotonicity method
of \cite{AFM}.
All inner products, gradients, divergences, and hypersurface measures
in Sections~4.1--4.2 use the metric~$g$. The main strategy is that using the conformal identities
of Proposition~\ref{prop:kato} and cutoff arguments across the
critical set, we establish the effective inequalities
$\Phi'(0)\le0$ and $\Phi(\infty)\le\Phi(0)$ in
Theorems~\ref{thm:firsteffective} and~\ref{thm:secondeffective}.
These yield the flux--curvature inequality by H\"older's inequality;
the equality case follows from the vanishing of the nonnegative
terms in the first divergence identity.

\subsection{First effective inequality: \texorpdfstring{$\Phi'(0)\le0$}{Phi'(0) <= 0}}

On the noncritical set, define
\begin{equation}\label{eq:X}
 X=u|\nabla\psi|^{p-k-2}
 \left[(n-k)|\nabla\psi|^2T_{k-1}(\KK)-kT_k(\KK)\right]\nabla_g\psi.
\end{equation}

\paragraph{Using  $ \Div_g J=0$ to compute $\Phi'$.}
On $\mathcal R$, the normal and mixed Newton blocks in \eqref{eq:blocknewton} give
$$
 \begin{aligned}
 (T_{k-1}(\KK)\nabla_g\psi)_\nu
 &=|\nabla\psi|^k\sigma_{k-1}(\HH),\\
 (T_{k-1}(\KK)\nabla_g\psi)^\top
 &=-|\nabla\psi|^{k-1}T_{k-2}(\HH)\nabla^\top |\nabla\psi|.
 \end{aligned}
$$
Multiplying by the factor $|\nabla\psi|^{p-k-1}$ in $J$, and then pairing
with $\nabla_g|\nabla\psi|=\nabla^\top |\nabla\psi|+\nu(|\nabla\psi|)\nu$, yields
\begin{equation}\label{eq:Jcomponents}
 \begin{aligned}
 J^\top&=-|\nabla\psi|^{p-2}T_{k-2}(\HH)\nabla^\top |\nabla\psi|,\\
 J_\nu&=|\nabla\psi|^{p-1}\sigma_{k-1}(\HH),\\
 \ip{J}{\nabla_g|\nabla\psi|}
 &=|\nabla\psi|^{p-2}\Bigl[|\nabla\psi|\sigma_{k-1}(\HH)\nu(|\nabla\psi|)
 -T_{k-2}(\HH)(\nabla^\top |\nabla\psi|,\nabla^\top |\nabla\psi|)\Bigr].
 \end{aligned}
\end{equation}
A subscript~$\nu$ here denotes a scalar normal component, not a
covariant coordinate. The product rule and $\Div_gJ=0$ imply
$\Div_g(|\nabla\psi|J)=\ip{J}{\nabla_g|\nabla\psi|}$. Substituting the level-frame
equation \eqref{eq:levelPDE} replaces the bracket in
\eqref{eq:Jcomponents} by an expression involving only $|\nabla\psi|$
and the level curvatures. Thus
\begin{equation}\label{eq:divqJ}
 \begin{aligned}
 \Div_g(|\nabla\psi|J)&=\ip{J}{\nabla_g|\nabla\psi|}\\
 &=\frac{|\nabla\psi|^p\sigma_{k-1}(\HH)}{p-1}
 \left[(n-k)|\nabla\psi|-k\frac{\sigma_k(\HH)}{\sigma_{k-1}(\HH)}\right].
 \end{aligned}
\end{equation}
In particular,
$$
 \int_{\Sigma_r}\ip{|\nabla\psi|J}{\nu}\dd S_g
 =\int_{\Sigma_r}|\nabla\psi|^p\sigma_{k-1}(\HH)\dd S_g=\Phi(r)
$$
whenever $r$ is a globally regular value.

For each globally regular value $s>0$, compactness of nearby
closed slabs gives a noncritical band around $\Sigma_s$.
The divergence theorem and coarea formula, followed by
\eqref{eq:divqJ}, therefore yield
\begin{equation}\label{eq:Phiderivative}
 \begin{aligned}
 (p-1)\Phi'(s)
 &=(p-1)\int_{\Sigma_s}
   \frac{\Div_g(|\nabla\psi|J)}{|\nabla\psi|}\dd S_g\\
 &=(n-k)\Phi(s)
   -k\int_{\Sigma_s}|\nabla\psi|^{p-1}
        \sigma_k(\HH)\dd S_g.
 \end{aligned}
\end{equation}
At $s=0$, the noncritical boundary collar in \textup{(P2)}
gives the same formula for the right derivative.
This calculation is local to the noncritical set and does
not use \textup{(P4)}.

\paragraph{The flux of $X$.}
The same Newton blocks, now with $j=k$, read
$$
 \begin{aligned}
 (T_k(\KK)\nabla_g\psi)_\nu&=|\nabla\psi|^{k+1}\sigma_k(\HH),\\
 (T_k(\KK)\nabla_g\psi)^\top&=-|\nabla\psi|^kT_{k-1}(\HH)\nabla^\top |\nabla\psi|.
 \end{aligned}
$$
Applying these and the corresponding $j=k-1$ blocks to
\eqref{eq:X} gives
\begin{equation}\label{eq:Xcomponents}
 \begin{aligned}
 X^\top&=u|\nabla\psi|^{p-2}
 \bigl[kT_{k-1}(\HH)-(n-k)|\nabla\psi|T_{k-2}(\HH)\bigr]\nabla^\top |\nabla\psi|,\\
 X_\nu&=u|\nabla\psi|^{p-1}
 \bigl[(n-k)|\nabla\psi|\sigma_{k-1}(\HH)-k\sigma_k(\HH)\bigr].
 \end{aligned}
\end{equation}
Since $\psi=s$ on $\Sigma_s$, equation
\eqref{eq:Phiderivative} now implies
\begin{equation}\label{eq:Xflux}
 \int_{\Sigma_s}\ip{X}{\nu}\dd S_g
 =(p-1)e^{-\frac{n-k-p+1}{p-1}s}\Phi'(s).
\end{equation}

\begin{lemma}%
\label{lem:X}
Assume \textup{(P1)--(P3)}. At every noncritical point,
\begin{align}
 \Div_gX
 &=u\Biggl\{
 k|\nabla\psi|^{p-1}\left[
  \frac{k(n-k-1)}{n-k}
    \frac{\sigma_k(\HH)^2}{\sigma_{k-1}(\HH)}
       -(k+1)\sigma_{k+1}(\HH)\right]\notag\\
 &\qquad +(p-1)k|\nabla\psi|^{p-3}
  \left[T_{k-1}(\HH)
    -\frac{\sigma_k(\HH)}{\sigma_{k-1}(\HH)}T_{k-2}(\HH)\right]
        (\nabla^\top |\nabla\psi|,\nabla^\top |\nabla\psi|)\notag\\
 &\qquad +\frac{|\nabla\psi|^{p-1}\sigma_{k-1}(\HH)}{n-k}
  \left[(n-k)|\nabla\psi|-k\frac{\sigma_k(\HH)}{\sigma_{k-1}(\HH)}\right]^2
 \Biggr\}\ge0,
 \label{eq:divX}\\
 \ip{X}{\nabla_g|\nabla\psi|}
 &=u\Biggl\{
 k|\nabla\psi|^{p-2}\left[T_{k-1}(\HH)
   -\frac{\sigma_k(\HH)}{\sigma_{k-1}(\HH)}T_{k-2}(\HH)\right]
       (\nabla^\top |\nabla\psi|,\nabla^\top |\nabla\psi|)\notag\\
 &\qquad+\frac{|\nabla\psi|^p\sigma_{k-1}(\HH)}{p-1}
  \left[(n-k)|\nabla\psi|-k\frac{\sigma_k(\HH)}{\sigma_{k-1}(\HH)}\right]^2
 \Biggr\}\ge0.
 \label{eq:contactX}
\end{align}
The first two terms in the divergence use precisely the scalar
and tensor curvature inequalities in \eqref{eq:BE}. The remaining
term is a square with a positive coefficient.
\end{lemma}

\begin{proof} %
We first write
$$
 \begin{aligned}
 X={}&(n-k)u|\nabla\psi|^{p-k}
             T_{k-1}(\KK)\nabla_g\psi -k u|\nabla\psi|^{p-k-2}T_k(\KK)\nabla_g\psi.
 \end{aligned}
$$
For the first summand, apply \eqref{eq:universaldiv} with
$\mu=\frac{n-k-p+1}{p-1}$, $\alpha=p-k$, and $j=k-1$.
The block formulas \eqref{eq:blockminor}--\eqref{eq:blocknewton}
give
\begin{align*}
 &u^{-1}\Div_g\!\left(
    u|\nabla\psi|^{p-k}T_{k-1}(\KK)\nabla_g\psi\right)\\
 &\quad=k|\nabla\psi|^p\sigma_k(\HH)
       -(n-k+\frac{n-k-p+1}{p-1})|\nabla\psi|^{p+1}\sigma_{k-1}(\HH)\\
 &\qquad+p |\nabla\psi|^{p-2}
   \left[|\nabla\psi|\sigma_{k-1}(\HH)\nu(|\nabla\psi|)
          -T_{k-2}(\HH)(\nabla^\top |\nabla\psi|,\nabla^\top |\nabla\psi|)\right].
\end{align*}
Here the two contributions to the coefficient $p$ come from the
principal-minor term and the derivative of the weight $|\nabla\psi|^{p-k}$.
For the second summand, use $\mu=\frac{n-k-p+1}{p-1}$, $\alpha=p-k-2$, and $j=k$:
\begin{align*}
 &u^{-1}\Div_g\!\left(
    u|\nabla\psi|^{p-k-2}T_k(\KK)\nabla_g\psi\right)\\
 &\quad=(k+1)|\nabla\psi|^{p-1}\sigma_{k+1}(\HH)
       -(n-k+\frac{n-k-p+1}{p-1}-1)|\nabla\psi|^p\sigma_k(\HH)\\
 &\qquad+(p-1)|\nabla\psi|^{p-3}
   \left[|\nabla\psi|\sigma_k(\HH)\nu(|\nabla\psi|)
          -T_{k-1}(\HH)(\nabla^\top |\nabla\psi|,\nabla^\top |\nabla\psi|)\right].
\end{align*}
This gives
\begin{align}
 u^{-1}|\nabla\psi|^{3-p}\Div_gX
 ={}&|\nabla\psi|^2\Bigl[
 k\bigl(2(n-k)+\frac{n-k-p+1}{p-1}-1\bigr)|\nabla\psi|\sigma_k(\HH)\notag\\
 &\qquad -(n-k)(n-k+\frac{n-k-p+1}{p-1})|\nabla\psi|^2\sigma_{k-1}(\HH)
             -k(k+1)\sigma_{k+1}(\HH)\Bigr]\notag\\
 &+(n-k)p |\nabla\psi|\left[
 |\nabla\psi|\sigma_{k-1}(\HH)\nu(|\nabla\psi|)
       -T_{k-2}(\HH)(\nabla^\top |\nabla\psi|,\nabla^\top |\nabla\psi|)\right]\notag\\
 &-k(p-1)\left[
 |\nabla\psi|\sigma_k(\HH)\nu(|\nabla\psi|)
       -T_{k-1}(\HH)(\nabla^\top |\nabla\psi|,\nabla^\top |\nabla\psi|)\right].
 \label{eq:divXexpanded}
\end{align}

We can now eliminate the normal derivative $\nu(|\nabla\psi|)$ using \eqref{eq:levelPDE}. In fact,
the two brackets containing $\nu(|\nabla\psi|)$ are related by
 \begin{align}
&|\nabla\psi|\sigma_k(\HH)\nu(|\nabla\psi|)
       -T_{k-1}(\HH)(\nabla^\top |\nabla\psi|,\nabla^\top |\nabla\psi|)\notag\\
 &\quad=\frac{\sigma_k(\HH)}{\sigma_{k-1}(\HH)}
    \left[|\nabla\psi|\sigma_{k-1}(\HH)\nu(|\nabla\psi|)
       -T_{k-2}(\HH)(\nabla^\top |\nabla\psi|,\nabla^\top |\nabla\psi|)\right]\notag\\
 &\qquad-\left[T_{k-1}(\HH)
       -\frac{\sigma_k(\HH)}{\sigma_{k-1}(\HH)}T_{k-2}(\HH)\right]
              (\nabla^\top |\nabla\psi|,\nabla^\top |\nabla\psi|).\label{eq:4.1}
 \end{align}
Thus \eqref{eq:levelPDE} says
 \begin{align}
 &(p-1)\left[|\nabla\psi|\sigma_{k-1}(\HH)\nu(|\nabla\psi|)
       -T_{k-2}(\HH)(\nabla^\top |\nabla\psi|,\nabla^\top |\nabla\psi|)\right]\notag\\
 &\qquad=|\nabla\psi|^2\sigma_{k-1}(\HH)
       \left[(n-k)|\nabla\psi|
          -k\frac{\sigma_k(\HH)}{\sigma_{k-1}(\HH)}\right].\label{eq:4.2}
 \end{align}
Substituting \eqref{eq:4.1} and \eqref{eq:4.2} into the last two lines of
\eqref{eq:divXexpanded} leaves the tensor term
$$
 (p-1)k\left[T_{k-1}(\HH)
       -\frac{\sigma_k(\HH)}{\sigma_{k-1}(\HH)}T_{k-2}(\HH)\right]
              (\nabla^\top |\nabla\psi|,\nabla^\top |\nabla\psi|)
$$
and the following scalar product, with no remaining $\nu(|\nabla\psi|)$:
$$
 \begin{aligned}
 &|\nabla\psi|^2\sigma_{k-1}(\HH)
   \left[\frac{(n-k)p}{p-1}|\nabla\psi|
         -k\frac{\sigma_k(\HH)}{\sigma_{k-1}(\HH)}\right]
   \left[(n-k)|\nabla\psi|
         -k\frac{\sigma_k(\HH)}{\sigma_{k-1}(\HH)}\right]\\
 &\quad=|\nabla\psi|^2\Biggl[
      \frac{(n-k)^2p}{p-1}|\nabla\psi|^2\sigma_{k-1}(\HH)
      -\left(\frac{(n-k)pk}{p-1}+k(n-k)\right)|\nabla\psi|\sigma_k(\HH)\\
 &\hspace{36mm}
      +k^2\frac{\sigma_k(\HH)^2}{\sigma_{k-1}(\HH)}\Biggr].
 \end{aligned}
$$
Using $\frac{n-k-p+1}{p-1}$, the total coefficient of
$|\nabla\psi|^2\sigma_{k-1}(\HH)$ in the curvature bracket becomes
$$
 \begin{aligned}
 &-(n-k)(n-k+\frac{n-k-p+1}{p-1})+\frac{(n-k)^2p}{p-1}=n-k.
 \end{aligned}
$$
The total coefficient of $|\nabla\psi|\sigma_k(\HH)$ becomes
$$
 \begin{aligned}
 &k\bigl(2(n-k)+\frac{n-k-p+1}{p-1}-1\bigr)
         -\frac{(n-k)pk}{p-1}-k(n-k)=-2k.\\
 \end{aligned}
$$
Consequently,
\begin{equation}\label{eq:divXpresquare}
 \begin{aligned}
 u^{-1}|\nabla\psi|^{3-p}\Div_gX
 ={}&|\nabla\psi|^2\Biggl[(n-k)|\nabla\psi|^2\sigma_{k-1}(\HH)-2k|\nabla\psi|\sigma_k(\HH)\\
 &\qquad+k^2\frac{\sigma_k(\HH)^2}{\sigma_{k-1}(\HH)}
                  -k(k+1)\sigma_{k+1}(\HH)\Biggr]\\
 &+(p-1)k\left[T_{k-1}(\HH)
       -\frac{\sigma_k(\HH)}{\sigma_{k-1}(\HH)}T_{k-2}(\HH)\right]
              (\nabla^\top |\nabla\psi|,\nabla^\top |\nabla\psi|).
 \end{aligned}
\end{equation}
The normal derivatives have been eliminated.

Next, we expand the square that appears in the statement:
$$
 \begin{aligned}
 &\frac{\sigma_{k-1}(\HH)}{n-k}
   \left[(n-k)|\nabla\psi|-k\frac{\sigma_k(\HH)}{\sigma_{k-1}(\HH)}\right]^2\\
 &\qquad=(n-k)|\nabla\psi|^2\sigma_{k-1}(\HH)-2k|\nabla\psi|\sigma_k(\HH)
       +\frac{k^2}{n-k}
             \frac{\sigma_k(\HH)^2}{\sigma_{k-1}(\HH)}.
 \end{aligned}
$$
Subtracting this from the curvature bracket in
\eqref{eq:divXpresquare} leaves
$$
 k\left[
   \frac{k(n-k-1)}{n-k}
       \frac{\sigma_k(\HH)^2}{\sigma_{k-1}(\HH)}
          -(k+1)\sigma_{k+1}(\HH)\right].
$$
This is $k$ times the nonnegative scalar expression in
\eqref{eq:BE}. Multiplying \eqref{eq:divXpresquare} by
$u|\nabla\psi|^{p-3}$ proves \eqref{eq:divX}.%

Finally, we will show $\ip{X}{\nabla_g|\nabla\psi|}_g\geq 0$. In fact, %
$$
 \begin{aligned}
 u^{-1}|\nabla\psi|^{2-p}\ip{X}{\nabla_g|\nabla\psi|}_g
 ={}&\left[kT_{k-1}(\HH)-(n-k)|\nabla\psi|T_{k-2}(\HH)\right]
                    (\nabla^\top |\nabla\psi|,\nabla^\top |\nabla\psi|)\\
 &+|\nabla\psi|\sigma_{k-1}(\HH)\nu(|\nabla\psi|)
   \left[(n-k)|\nabla\psi|-k\frac{\sigma_k(\HH)}{\sigma_{k-1}(\HH)}\right].
 \end{aligned}
$$
Use \eqref{eq:levelPDE} in the form
$$
 \begin{aligned}
 |\nabla\psi|\sigma_{k-1}(\HH)\nu(|\nabla\psi|)
 ={}&T_{k-2}(\HH)(\nabla^\top |\nabla\psi|,\nabla^\top |\nabla\psi|)\\
 &+\frac{|\nabla\psi|^2\sigma_{k-1}(\HH)}{p-1}
   \left[(n-k)|\nabla\psi|-k\frac{\sigma_k(\HH)}{\sigma_{k-1}(\HH)}\right].
 \end{aligned}
$$
After substitution, the coefficient of
$T_{k-2}(\HH)(\nabla^\top |\nabla\psi|,\nabla^\top |\nabla\psi|)$ is
$$
 \left[(n-k)|\nabla\psi|-k\frac{\sigma_k(\HH)}{\sigma_{k-1}(\HH)}\right]
 -(n-k)|\nabla\psi|
 =-k\frac{\sigma_k(\HH)}{\sigma_{k-1}(\HH)}.
$$
Combining this term with
$kT_{k-1}(\HH)(\nabla^\top |\nabla\psi|,\nabla^\top |\nabla\psi|)$ gives
$$
 \begin{aligned}
 &u^{-1}|\nabla\psi|^{2-p}\ip{X}{\nabla_g|\nabla\psi|}\\
 &\quad=k\left[T_{k-1}(\HH)
       -\frac{\sigma_k(\HH)}{\sigma_{k-1}(\HH)}T_{k-2}(\HH)\right]
              (\nabla^\top |\nabla\psi|,\nabla^\top |\nabla\psi|)\\
 &\qquad+\frac{|\nabla\psi|^2\sigma_{k-1}(\HH)}{p-1}
   \left[(n-k)|\nabla\psi|-k\frac{\sigma_k(\HH)}{\sigma_{k-1}(\HH)}\right]^2.
 \end{aligned}
$$
Multiplying by $u|\nabla\psi|^{p-2}$ proves \eqref{eq:contactX}.
Finally, $|\nabla\psi|>0$, $p-1>0$, $n-k>0$, and
$\sigma_{k-1}(\HH)>0$. Together with the two inequalities in
\eqref{eq:BE}, we know every term in
\eqref{eq:divX} and \eqref{eq:contactX} are nonnegative.
\end{proof}

\begin{remark}    
For $k=1$, we have $\sigma_0(\HH)=1$, $T_{-1}(\HH)=0$,
$T_0(\HH)=I$, and the equation is $\Delta_p^g\psi=0$.
Thus \eqref{eq:levelPDE} reads
$$
 (n-1)|\nabla\psi|-\sigma_1(\HH)=\frac{p-1}{|\nabla\psi|}\nu(|\nabla\psi|).
$$
Thus
\begin{align*}
 X&=u|\nabla\psi|^{p-2}
       \left[\nabla^\top |\nabla\psi|+(p-1)\nu(|\nabla\psi|)\nu\right],\\
 \Div_gX
 &=u|\nabla\psi|^{p-3}
   \left[|\nabla\psi|^2|h^\circ|_g^2+(p-1)\bigl|\nabla^\top|\nabla\psi|\bigr|_g^2
                    +\frac{(p-1)^2}{n-1}\nu(|\nabla\psi|)^2\right].
\end{align*}
The higher-order identity therefore reduces to the trace-free,
tangential, and normal decomposition used in the $p$-Laplace case.
\end{remark}

\begin{theorem}%
\label{thm:firsteffective}
Assume \textup{(P1)--(P3)}. For two globally regular values
$0\le s<t$,
\begin{equation}\label{eq:firsteffective}
 e^{-\frac{n-k-p+1}{p-1}s}\Phi'(s)\le e^{-\frac{n-k-p+1}{p-1}t}\Phi'(t).
\end{equation}
Critical points and critical values are allowed between the two levels.
If, in addition, all sufficiently large levels are regular and $\Phi$
is bounded at infinity, then
\begin{equation}\label{eq:firstendpoint}
 \Phi'(s)\le0
 \quad\text{at every globally regular value }s\ge0.
 \qquad\text{In particular, }\Phi'(0)\le0.
\end{equation}
The additional conditions follow from \textup{(P5)}.
\end{theorem}

\begin{proof}
Fix globally regular
values $0\le s<t$. The slab
$\mathcal A_{s,t}=\{s<\psi<t\}$ has compact closure, and
\[
 m_{s,t}:=\min_{\Sigma_s\cup\Sigma_t}|\nabla\psi|>0.
\]
When $s=0$, we interpret $\Sigma_0=\partial\Omega$.

Choose $0<2\varepsilon<m_{s,t}$ and a smooth nondecreasing
function $\chi_\varepsilon:[0,\infty)\to[0,1]$ such that
\[
 \chi_\varepsilon=0\quad\text{on }[0,\varepsilon],
 \qquad
 \chi_\varepsilon=1\quad\text{on }[2\varepsilon,\infty).
\]
Define $X_\varepsilon=\chi_\varepsilon(|\nabla\psi|)X$
on $\mathcal R$ and extend it by zero across the critical set.
By continuity of $|\nabla\psi|$, this field vanishes in a
neighborhood of every critical point. Hence $X_\varepsilon$
is $C^1$ on the closed slab and agrees with $X$ on its boundary.

By \eqref{eq:divX} and \eqref{eq:contactX}, on $\mathcal R$,
\begin{equation}\label{eq:Xcutoffdiv}
 \Div_gX_\varepsilon
 =\chi_\varepsilon(|\nabla\psi|)\Div_gX
 +\chi_\varepsilon'(|\nabla\psi|)
     \ip{X}{\nabla_g|\nabla\psi|}
 \ge0.
\end{equation}
The divergence is zero near the critical set. The divergence
theorem and \eqref{eq:Xflux} therefore give
\begin{align*}
 0
 &\le\int_{\mathcal A_{s,t}}\Div_gX_\varepsilon\dd V_g\\
 &=\int_{\Sigma_t}\ip{X}{\nu}\dd S_g
   -\int_{\Sigma_s}\ip{X}{\nu}\dd S_g\\
 &=(p-1)\bigl[e^{-\alpha_p t}\Phi'(t)
                  -e^{-\alpha_p s}\Phi'(s)\bigr],
\end{align*}
which proves \eqref{eq:firsteffective}. This argument uses
a fixed sufficiently small $\varepsilon$ and does not
require \textup{(P4)}.

Now assume that all sufficiently large levels are regular
and $\Phi$ is bounded at infinity. If $\Phi'(s)>0$ at a
globally regular value $s$, choose $T>s$ beyond which all
levels are regular. By \eqref{eq:firsteffective},
\[
 \Phi'(r)\ge e^{\alpha_p(r-s)}\Phi'(s),
 \qquad r\ge T.
\]
Since $\Phi$ is $C^1$ on the regular tail, integration yields
\[
 \Phi(t)-\Phi(T)
 \ge \frac{e^{-\alpha_p s}\Phi'(s)}{\alpha_p}
       \bigl(e^{\alpha_p t}-e^{\alpha_p T}\bigr),
 \qquad t>T,
\]
contradicting boundedness as $t\to\infty$. Thus
$\Phi'(s)\le0$ at every globally regular value, including
the right derivative at $s=0$ supplied by the boundary
collar in \textup{(P2)}. Under \textup{(P5)}, the required
regularity and boundedness at infinity follow from the
previous asymptotic analysis.
\end{proof}

\subsection{Second effective inequality: \texorpdfstring{$\Phi(\infty)\le\Phi(0)$}{Phi(infinity) <= Phi(0)}}

For $0<\lambda<1$, define on $\mathcal R$
\begin{equation}\label{eq:Y}
 Y_\lambda=(u^{-1}-\lambda)X-(n-k-p+1)|\nabla\psi|J,
 \qquad
 \eta_\lambda=(u^{-1}-\lambda)^{-1}>0.
\end{equation}
The correction term $(n-k-p+1)|\nabla\psi|J$ is chosen to cancel the derivative of $u^{-1}$ in
the divergence. The following lemma shows that $Y_\lambda$
has nonnegative divergence and that nondecreasing cutoffs
in $\eta_\lambda|\nabla\psi|$ preserve this property.

\begin{lemma}%
\label{lem:Y}
Under \textup{(P1)--(P3)}, the following identities hold on the
regular set for $0<\lambda<1$:
\begin{equation}\label{eq:divY}
 \Div_gY_\lambda=(u^{-1}-\lambda)\Div_gX\ge0.
\end{equation}
\begin{equation}\label{eq:contactY}
 \begin{aligned}
 \ip{Y_\lambda}{\nabla_g(\eta_\lambda |\nabla\psi|)}
 &\quad=k u|\nabla\psi|^{p-2}
 \left[T_{k-1}(\HH)
 -\frac{\sigma_k(\HH)}{\sigma_{k-1}(\HH)}T_{k-2}(\HH)\right]
 (\nabla^\top |\nabla\psi|,\nabla^\top |\nabla\psi|)\\
 &\qquad+\frac{u|\nabla\psi|^p\sigma_{k-1}(\HH)}{p-1}
 \left[
 \begin{aligned}
 &(n-k)|\nabla\psi|-k\frac{\sigma_k(\HH)}{\sigma_{k-1}(\HH)}\\
 &\quad-(n-k-p+1)u^{-1}\eta_\lambda |\nabla\psi|
 \end{aligned}
 \right]^2\ge0.
 \end{aligned}
\end{equation}
\end{lemma}

\begin{proof}
Set
\[
 q=|\nabla\psi|,\qquad
 B=(n-k)q-k\frac{\sigma_k(\HH)}{\sigma_{k-1}(\HH)}.
\]
Since $\nabla_g\psi=q\nu$, equations \eqref{eq:Xcomponents}
and \eqref{eq:divqJ} give
\begin{equation}\label{eq:XYcancellation}
 u^{-1}\ip{X}{\nabla_g\psi}
 =q^p\sigma_{k-1}(\HH)B
 =(p-1)\Div_g(qJ).
\end{equation}
Using
\[
 \nabla_g(u^{-1})
 =\frac{n-k-p+1}{p-1}u^{-1}\nabla_g\psi,
\]
the product rule yields
\begin{align*}
 \Div_gY_\lambda
 &=(u^{-1}-\lambda)\Div_gX
   +\frac{n-k-p+1}{p-1}u^{-1}\ip{X}{\nabla_g\psi}
   -(n-k-p+1)\Div_g(qJ)\\
 &=(u^{-1}-\lambda)\Div_gX\ge0.
\end{align*}
This proves \eqref{eq:divY}, since $u^{-1}-\lambda>0$
and $\Div_gX\ge0$ by \eqref{eq:divX}.

Next,
\[
 \nabla_g\eta_\lambda
 =-\frac{n-k-p+1}{p-1}
   u^{-1}\eta_\lambda^2\nabla_g\psi.
\]
Since $Y_\lambda=\eta_\lambda^{-1}X-(n-k-p+1)qJ$,
we obtain
\begin{align*}
 \ip{Y_\lambda}{\nabla_g(\eta_\lambda q)}
 ={}&\ip{X}{\nabla_gq}\\
 &-\frac{n-k-p+1}{p-1}
   u^{-1}\eta_\lambda q\ip{X}{\nabla_g\psi}\\
 &-(n-k-p+1)\eta_\lambda\ip{qJ}{\nabla_gq}\\
 &+\frac{(n-k-p+1)^2}{p-1}
   u^{-1}\eta_\lambda^2q\ip{qJ}{\nabla_g\psi}.
\end{align*}
Equations \eqref{eq:Jcomponents} and \eqref{eq:divqJ} imply
\[
 \begin{aligned}
 \ip{qJ}{\nabla_gq}
 &=\frac{q^{p+1}\sigma_{k-1}(\HH)}{p-1}B,\\
 \ip{qJ}{\nabla_g\psi}
 &=q^{p+1}\sigma_{k-1}(\HH).
 \end{aligned}
\]
Together with \eqref{eq:XYcancellation}, these identities show
that the two mixed terms coincide. Thus
\begin{align*}
 \ip{Y_\lambda}{\nabla_g(\eta_\lambda q)}
 ={}&\ip{X}{\nabla_gq}
 -\frac{2(n-k-p+1)}{p-1}
   \eta_\lambda q^{p+1}\sigma_{k-1}(\HH)B\\
 &+\frac{(n-k-p+1)^2}{p-1}
   u^{-1}\eta_\lambda^2q^{p+2}\sigma_{k-1}(\HH).
\end{align*}
Substituting \eqref{eq:contactX} and completing the square gives
\[
 \begin{aligned}
 \ip{Y_\lambda}{\nabla_g(\eta_\lambda q)}
 ={}&kuq^{p-2}
 \left[
 T_{k-1}(\HH)
 -\frac{\sigma_k(\HH)}{\sigma_{k-1}(\HH)}T_{k-2}(\HH)
 \right](\nabla^\top q,\nabla^\top q)\\
 &+\frac{uq^p\sigma_{k-1}(\HH)}{p-1}
 \left[B-(n-k-p+1)u^{-1}\eta_\lambda q\right]^2.
 \end{aligned}
\]
This is \eqref{eq:contactY}. The tangential term is nonnegative
by \eqref{eq:BE}, and the square has a positive coefficient
because $u,q,\sigma_{k-1}(\HH)>0$ and $p>1$.
\end{proof}

\begin{theorem}%
\label{thm:secondeffective}
Assume \textup{(P1)--(P3)}. For any two globally regular values
$0\le s<t$ and every $0<\lambda<1$,
\begin{equation}\label{eq:secondeffective}
 (1-\lambda e^{-\frac{n-k-p+1}{p-1}s})\Phi'(s)-\frac{n-k-p+1}{p-1}\Phi(s)
 \le
 (1-\lambda e^{-\frac{n-k-p+1}{p-1}t})\Phi'(t)-\frac{n-k-p+1}{p-1}\Phi(t).
\end{equation}
If \textup{(P5)} is also satisfied, then
\begin{equation}\label{eq:secondendpoint}
 \Phi(\infty)\le\Phi(0).
\end{equation}
Neither assertion requires \textup{(P4)}.
\end{theorem}

\begin{proof}
We first construct a cutoff function.
Fix the two regular endpoints $s<t$ and a number
$0<\lambda<1$. Throughout their closed slab,
$$
 0<\frac{1}{e^{\frac{n-k-p+1}{p-1}t}-\lambda}\le\eta_\lambda
 \le\frac{1}{e^{\frac{n-k-p+1}{p-1}s}-\lambda}<\infty.
$$
The function $\eta_\lambda |\nabla\psi|$ is therefore continuous, is zero
at every critical point, and is positive on both regular boundary
levels. Let
$m_{s,t,\lambda}=\min_{\Sigma_s\cup\Sigma_t}(\eta_\lambda |\nabla\psi|)>0$.
Choose the nondecreasing cutoff from the preceding theorem with
$0<2\varepsilon<m_{s,t,\lambda}$. Define the following field
on the regular set and extend it by zero at critical points:
$$
 Y_{\lambda,\varepsilon}
 =\chi_\varepsilon(\eta_\lambda |\nabla\psi|)Y_\lambda.
$$

As in the first cutoff argument, continuity gives a neighborhood of
each critical point on which $\eta_\lambda |\nabla\psi|<\varepsilon$.
The cut field is identically zero in that neighborhood; away from
the critical set it is $C^1$. Thus it is $C^1$ on the entire
closed slab, regardless of whether the uncut~$Y_\lambda$ extends
to its critical points.

The product rule and Lemma~\ref{lem:Y} give
\begin{align*}
 \Div_gY_{\lambda,\varepsilon}
 ={}&\chi_\varepsilon(\eta_\lambda |\nabla\psi|)\Div_gY_\lambda\\
 &+\chi_\varepsilon'(\eta_\lambda |\nabla\psi|)
   \ip{Y_\lambda}{\nabla_g(\eta_\lambda |\nabla\psi|)}\ge0.
\end{align*}
Both signs hold on every noncritical point of the slab, including
points whose level values are critical elsewhere. This is where the
pointwise scope of \textup{(P3)} is needed. Near the critical set
the cut field has zero divergence.

For a regular value $r$, \eqref{eq:Xflux} and the flux of~$|\nabla\psi|J$
give
\begin{align}
 \int_{\Sigma_r}\ip{Y_\lambda}{\nu}\dd S_g
 &=(e^{\frac{n-k-p+1}{p-1}r}-\lambda)(p-1)e^{-\frac{n-k-p+1}{p-1}r}\Phi'(r)-(n-k-p+1)\Phi(r)\notag\\
 &=(p-1)\left[(1-\lambda e^{-\frac{n-k-p+1}{p-1}r})\Phi'(r)-\frac{n-k-p+1}{p-1}\Phi(r)\right].
 \label{eq:Yflux}
\end{align}
Our fixed choice of $\varepsilon$ makes the cutoff equal to one
on both endpoints, so these are also the boundary fluxes of
$Y_{\lambda,\varepsilon}$. Applying the divergence theorem,
with $+\nu$ at $\Sigma_t$ and $-\nu$ at $\Sigma_s$, gives
\begin{align*}
 0&\le\int_{\mathcal A_{s,t}}
                   \Div_gY_{\lambda,\varepsilon}\dd V_g\\
 &=(p-1)\bigl[(1-\lambda e^{-\frac{n-k-p+1}{p-1}t})\Phi'(t)-\frac{n-k-p+1}{p-1}\Phi(t)\\
 &\hspace{37mm}-(1-\lambda e^{-\frac{n-k-p+1}{p-1}s})\Phi'(s)+\frac{n-k-p+1}{p-1}\Phi(s)\bigr].
\end{align*}
After division by $p-1>0$, this is
\eqref{eq:secondeffective}. The cutoff parameter has disappeared
from this scalar inequality; no limit of cut vector fields is taken.

Under \textup{(P5)}, every sufficiently large $t$ is regular,
$\Phi(t)\to\Phi(\infty)<\infty$, and
Theorem~\ref{thm:firsteffective} gives $\Phi'(t)\le0$.
Take $s=0$ in \eqref{eq:secondeffective}. Rearranging the terms,
$$
 \frac{n-k-p+1}{p-1}[\Phi(t)-\Phi(0)]
 \le(1-\lambda e^{-\frac{n-k-p+1}{p-1}t})\Phi'(t)-(1-\lambda)\Phi'(0).
$$
Keep this regular $t$ fixed. All quantities except $\lambda$
are finite scalar values. Letting $\lambda\uparrow1$ now gives
$$
 \frac{n-k-p+1}{p-1}[\Phi(t)-\Phi(0)]\le(1-e^{-\frac{n-k-p+1}{p-1}t})\Phi'(t)\le0.
$$
Because $\frac{n-k-p+1}{p-1}>0$, we have proved $\Phi(t)\le\Phi(0)$ on the
regular tail. Finally take $t\to\infty$ to obtain
\eqref{eq:secondendpoint}.

The cutoff argument is applied for each fixed $\lambda<1$, with
$\varepsilon=\varepsilon(s,t,\lambda)$, before passing to
$\lambda\uparrow1$ in the scalar inequality. This order avoids
the singular field at $\lambda=1$ near $\psi=0$ and requires
no uniform estimate for the cutoff fields. The subsequent limit
uses the sign of $\Phi'(t)$, rather than convergence of that derivative.
\end{proof}

\paragraph{Comparison between regular levels.}
On a single slab $\mathcal A_{s,t}$, the positivity condition
needed for the second field is only
$u^{-1}-\lambda>0$ there. Thus its proof also works for
$0<\lambda<e^{\frac{n-k-p+1}{p-1}s}$, since $\psi\ge s$. The field need not
be defined below that slab. The bound
$$
 (e^{\frac{n-k-p+1}{p-1}t}-\lambda)^{-1}\le\eta_\lambda
 \le(e^{\frac{n-k-p+1}{p-1}s}-\lambda)^{-1}
$$
remains valid for each such fixed~$\lambda$. The same
construction, with its cutoff chosen after fixing $\lambda$,
therefore proves \eqref{eq:secondeffective} in this enlarged range.
Rearrange it as
$$
 \frac{n-k-p+1}{p-1}[\Phi(t)-\Phi(s)]
 \le(1-\lambda e^{-\frac{n-k-p+1}{p-1}t})\Phi'(t)
    -(1-\lambda e^{-\frac{n-k-p+1}{p-1}s})\Phi'(s).
$$
Now let $\lambda\uparrow e^{\frac{n-k-p+1}{p-1}s}$ in the scalar inequality.
The coefficient of $\Phi'(s)$ tends to zero, whereas the
coefficient of $\Phi'(t)$ tends to
$1-e^{-\frac{n-k-p+1}{p-1}(t-s)}$. Hence, under \textup{(P5)},
\begin{equation}\label{eq:regularmonotonicity}
 \frac{n-k-p+1}{p-1}[\Phi(t)-\Phi(s)]
 \le\left(1-e^{-\frac{n-k-p+1}{p-1}(t-s)}\right)\Phi'(t)\le0.
\end{equation}
The last sign is Theorem~\ref{thm:firsteffective}. This proves
$\Phi(t)\le\Phi(s)$ for any two globally regular endpoints,
even if critical values lie between them.

This conclusion compares geometric integrals that are already defined
on regular levels. It neither proves that globally regular values are
dense nor identifies a monotone extension with the geometric level
integral at every critical value. We retain the term \emph{effective
inequalities} to distinguish these proved statements from an additional
all-level formulation. In particular, no argument above establishes
weak conservation of the original Newton flux across the critical set:
that remains the separate hypothesis \textup{(P4)} needed in the next
subsection. The two endpoint inequalities
$\Phi'(0)\le0$ and $\Phi(\infty)\le\Phi(0)$ already suffice
for the geometric application.

\subsection{The boundary inequality and rigidity}
In this subsection,  Assumptions~\ref{ass:potential} and~\ref{ass:asymptotic} are imposed, including
the weak Euclidean flux condition \textup{(P4)}.
The fixed-parameter estimate uses the boundary curvature conditions
through the positivity of $H_{k-1}$ and $H_k$. It does not require
convexity of the interior level hypersurfaces.

\begin{mainresult}%
Let $n\ge3$, $1\le k\le n-2$, and $1<p<n-k$.
Let $\Omega\subset\RR^n$ be a smooth bounded strictly star-shaped
domain with strictly $k$-convex boundary,
and let $u$ satisfy Assumptions~\ref{ass:potential}
and~\ref{ass:asymptotic}. Then its normalized flux satisfies
\begin{equation}\label{eq:fixedinequality}
 \CC(\Omega;u)^{\frac{n-k-p}{n-k-p+1}}
 \le \avgint_{\partial\Omega}
                \frac{H_k^p}{H_{k-1}^{p-1}}\dd S_\delta.
\end{equation}
Equality holds if and only if $\Omega$ is a ball. If
$\Omega=B_R(x_0)$, the potential is necessarily
$u(x)=(R/|x-x_0|)^{\alpha_p}$.
\end{mainresult}

\begin{proof}(Inequality case)
 Recall the normalized
curvatures and flux,
$$
 H_j=\binom{n-1}{j}^{-1}\sigma_j(h_\delta),\qquad
 \CC(\Omega)=\avgint_{\partial\Omega}
       \left(\frac{(p-1)|Du|}{n-k-p+1}\right)^{p-1}H_{k-1}\dd S_\delta.
$$
The symbol $\avgint$ means division by the fixed number
$\omega_{n-1}$. All applications of H\"older below use the
same normalized surface measure
$\omega_{n-1}^{-1}\dd S_\delta$.
Lemma~\ref{lem:flux} identifies the conserved flux with the
far-field coefficient, $\CC=\mathsf a^{p-1}$.
Recalling $\widehat\Phi=\Phi/[\omega_{n-1}\binom{n-1}{k-1}]$,
the limit \eqref{eq:Phiinfinity} and the second effective inequality give
\begin{equation}\label{eq:fluxlowerbound}
 \CC(\Omega)^{\frac{n-k-p}{n-k-p+1}}
 =\widehat\Phi(\infty)\le\widehat\Phi(0).
\end{equation}

On $\partial\Omega$, the boundary value $u=1$ gives
$\psi=0$, $g=\delta$, $|\nabla\psi|=\frac{(p-1)|Du|}{n-k-p+1}$, and $\HH=h_\delta$.
Thus all boundary integrals below use the Euclidean metric, and
$$
 \widehat\Phi(0)
 =\avgint_{\partial\Omega}|\nabla\psi|^pH_{k-1}\dd S_\delta.
$$
We will bound this integral using the first effective inequality,
then remove $|\nabla\psi|$ by H\"older's inequality.

To derive the boundary derivative inequality, we divide \eqref{eq:Phiderivative} at $s=0$ by
$\omega_{n-1}\binom{n-1}{k-1}$. The curvature term becomes
$$
 \begin{aligned}
 &\frac{k}{\omega_{n-1}\binom{n-1}{k-1}}
 \int_{\partial\Omega}|\nabla\psi|^{p-1}\sigma_k(h_\delta)\dd S_\delta\\
 &\qquad=(n-k)\avgint_{\partial\Omega}|\nabla\psi|^{p-1}H_k\dd S_\delta,
 \end{aligned}
$$
because
$k\binom{n-1}{k}=(n-k)\binom{n-1}{k-1}$. Hence
$$
 (p-1)\widehat\Phi'(0)
 =(n-k)\left[
 \widehat\Phi(0)-\avgint_{\partial\Omega}|\nabla\psi|^{p-1}H_k\dd S_\delta
 \right].
$$
Since $p-1>0$ and $n-k>0$, the inequality
$\widehat\Phi'(0)\le0$ gives
\begin{equation}\label{ineqn:4.3}
 \widehat\Phi(0)
 \le\avgint_{\partial\Omega}|\nabla\psi|^{p-1}H_k\dd S_\delta.
\end{equation}
By H\"older's inequality 
$$
 \avgint_{\partial\Omega}|\nabla\psi|^{p-1}H_k\dd S_\delta
 \le\widehat\Phi(0)^{(p-1)/p}
 \left(\avgint_{\partial\Omega}
       \frac{H_k^p}{H_{k-1}^{p-1}}\dd S_\delta\right)^{1/p}.
$$
Here we have used
$\widehat\Phi(0)>0$, and that
$H_{k-1}>0$, $H_k>0$, and $|\nabla\psi|>0$.

Combining this with \eqref{ineqn:4.3} gives
$$
 \widehat\Phi(0)
 \le\widehat\Phi(0)^{(p-1)/p}
 \left(\avgint_{\partial\Omega}
       \frac{H_k^p}{H_{k-1}^{p-1}}\dd S_\delta\right)^{1/p}.
$$
Thus
$$
 \widehat\Phi(0)
 \le\avgint_{\partial\Omega}
       \frac{H_k^p}{H_{k-1}^{p-1}}\dd S_\delta.
$$
Together with \eqref{eq:fluxlowerbound}, this proves the inequality \eqref{eq:fixedinequality}.

\end{proof}

\begin{proof}(Equality case)
For the explicit radial potential of a ball of radius~$R$,
\eqref{eq:ball} gives $\CC(B_R)=R^{n-k-p+1}$, so the left
side of \eqref{eq:fixedinequality} is $R^{n-k-p}$.
On its boundary $H_j=R^{-j}$ and
$\avgint_{\partial B_R}1\dd S_\delta=R^{n-1}$, so the right
side is also
$$
 R^{n-1}R^{-kp+(k-1)(p-1)}=R^{n-k-p}.
$$

\emph{Uniqueness of the potential on a ball.}
We show that every potential in the stated class on a ball is the
radial one. After translating the center to the origin, set
$$
 v(x)=\left(\frac R{|x|}\right)^{\alpha_p},
 \qquad \alpha_p=\frac{n-k-p+1}{p-1}.
$$
All derivatives, normals, and matrix inequalities in this uniqueness
argument are Euclidean. We use $P$ from
Section~\ref{sec:weightedhessian}; it depends on the gradient direction.
The argument uses \textup{(P1)}, \textup{(P2)}, \textup{(P4)}, and
\textup{(P5)}.

For $0<\gamma<1$, the gradients of $v^\gamma$ and $v$ have the same
direction and never vanish. In a radial orthonormal frame,
$$
 -P(Dv)D^2(v^\gamma)P(Dv)
 =\frac{\gamma\alpha_p v^\gamma}{|x|^2}
   \operatorname{diag}\left(
      1,\ldots,1,-\frac{p-1}{k}(\gamma\alpha_p+1)\right).
$$
For $1\le j\le k$, its $j$-th elementary symmetric polynomial is
$$
 \left(\frac{\gamma\alpha_p v^\gamma}{|x|^2}\right)^j
 \left[\binom{n-1}{j}
 -\frac{p-1}{k}(\gamma\alpha_p+1)\binom{n-1}{j-1}\right]>0,
$$
since
$$
 \frac{p-1}{k}(\gamma\alpha_p+1)
 <\frac{n-k}{k}\le\frac{n-j}{j}.
$$
Thus the weighted Hessian of $v^\gamma$ lies in $\Gamma_k$.
If $u-v^\gamma$ were positive somewhere, its zero boundary values
and decay at infinity would give a positive maximum at a finite
interior point. There $Du=D(v^\gamma)\ne0$, so \textup{(P2)} gives
the classical equation. The Hessian ordering yields
$$
 -P(Du)D^2uP(Du)
 \ge -P(Dv)D^2(v^\gamma)P(Dv)\in\Gamma_k.
$$
Addition of a positive semidefinite matrix preserves $\Gamma_k$,
contradicting $\sigma_k(-PD^2uP)=0$. Hence $u\le v^\gamma$;
letting $\gamma\uparrow1$ gives $u\le v$ throughout the exterior.
This comparison requires no branch assumption on $u$.

On $\partial B_R$, both functions equal one. Differentiating their
order along the outward normal of the ball gives
$$
 |Du|\ge |Dv|=\frac{\alpha_p}{R}.
$$
Here the tangential derivatives vanish and
$\partial_{\nu_\delta}u=-|Du|$, since $u<1$ in the exterior and
the boundary gradient is nonzero. At infinity,
$u\le v$ and \textup{(P5)} give $\mathsf a\le R^{\alpha_p}$.
Lemma~\ref{lem:flux}, using \textup{(P4)} and \textup{(P5)},
identifies the same flux at the boundary and at infinity. Since
$H_{k-1}=R^{-(k-1)}$ on the sphere,
$$
 \begin{aligned}
 R^{n-k-p+1}
 &\le \avgint_{\partial B_R}
       \left(\frac{|Du|}{\alpha_p}\right)^{p-1}
       R^{-(k-1)}\dd S_\delta\\
 &=\mathsf a^{p-1}\le R^{n-k-p+1}.
 \end{aligned}
$$
All inequalities are equalities. The boundary integrand is continuous
and pointwise no smaller than its radial value, so
$|Du|=\alpha_p/R$ everywhere on $\partial B_R$.
In particular $Du=Dv$ there. This conclusion uses only the flux
identity, not the inequality being proved.

Tangential differentiation of $Du=Dv$ on the sphere shows that all
Hessian components involving a tangent vector agree. The mixed
components vanish, and the tangential block of $-D^2u$ is
$(|Du|/R)I$. The classical equation extends to the boundary by
\textup{(P2)} and gives
$$
 0=\binom{n-1}{k}\left(\frac{|Du|}{R}\right)^k
 -\frac{p-1}{k}D^2u(\nu_\delta,\nu_\delta)
       \binom{n-1}{k-1}\left(\frac{|Du|}{R}\right)^{k-1}.
$$
Consequently,
$$
 D^2u(\nu_\delta,\nu_\delta)
 =\frac{n-k}{p-1}\frac{|Du|}{R}
 =D^2v(\nu_\delta,\nu_\delta).
$$
Thus the functions, their gradients, and their Hessians agree on
the entire sphere.

At the radial potential,
$$
 -P(Dv)D^2vP(Dv)
 =\frac{|Dv|}{|x|}
   \operatorname{diag}\left(1,\ldots,1,-\frac{n-k}{k}\right).
$$
The normal and tangential eigenvalues of its Newton tensor
$T_{k-1}$ are, respectively,
$$
 \binom{n-1}{k-1}\left(\frac{|Dv|}{|x|}\right)^{k-1},
 \qquad
 \frac1k\binom{n-2}{k-1}
       \left(\frac{|Dv|}{|x|}\right)^{k-1}.
$$
Both are positive, including when $k=1$. The derivative with respect
to the Hessian of the normalized equation
$-\sigma_k(-P(Du)D^2uP(Du))=0$ is $P T_{k-1}P$.
It is therefore positive definite near the radial gradient and
Hessian at every finite point.

The matching boundary derivatives and \textup{(P2)} allow us to
choose a thin spherical collar in which all line segments joining
the gradients and Hessians of $u$ and $v$ remain in this elliptic
neighborhood, with nonzero gradients. Integrating the derivatives
of the normalized operator along these segments shows that $u-v$
satisfies a linear uniformly elliptic equation in the collar.
Its first-order coefficients are continuous and bounded, and its
zeroth-order coefficient is zero because the operator does not
depend on the function value. These local bounds may depend on
the fixed $u,R,n,k,p$.

We have $u-v\le0$, with zero value and zero normal derivative on
the inner boundary of the collar. If the difference were not
identically zero, the strong maximum principle would make it
strictly negative in the collar. The Hopf boundary point lemma,
applied using a ball contained in the collar and tangent to its
inner boundary, would then give a nonzero normal derivative,
a contradiction. Here we use precisely the linear results
\cite[Lemma~3.4 and Theorem~3.5]{GilbargTrudinger2001}.
It follows that $u=v$ in the collar.

To extend this equality, consider the nonempty open set of exterior
points near which $u=v$ identically. At any of its limit points in
the exterior, \textup{(P1)} gives $Du=Dv\ne0$. Thus \textup{(P2)}
provides a $C^3$ neighborhood; taking limits from the equality set
also gives $D^2u=D^2v$ at that point. On a sufficiently small ball
around it, the same linearization is uniformly elliptic. Since
$u-v\le0$ attains the value zero there, the strong maximum principle
gives equality throughout that small ball. The equality set is
therefore also relatively closed. Connectedness of the exterior
implies $u=v$ everywhere. Each possible endpoint of this extension
is noncritical, so no regularity across an unknown critical set is
being assumed. Every potential in the stated class on a ball is
thus the radial potential and attains equality.

\emph{Rigidity for a general domain.}
We return to the conformal notation used above.
Conversely, equality in \eqref{eq:fixedinequality} forces equality
in both bounds for $\widehat\Phi(0)$, and therefore
$$
 \widehat\Phi(0)
 =\avgint_{\partial\Omega}
       \frac{H_k^p}{H_{k-1}^{p-1}}\dd S_\delta.
$$
Substituting this equality into the H\"older bound gives
$$
 \widehat\Phi(0)
 \le\avgint_{\partial\Omega}|\nabla\psi|^{p-1}H_k\dd S_\delta
 \le\widehat\Phi(0)^{(p-1)/p}\widehat\Phi(0)^{1/p}
 =\widehat\Phi(0).
$$
Thus the middle integral equals $\widehat\Phi(0)$. The boundary
derivative formula now yields $\Phi'(0)=0$.

For every globally regular~$t$,
Theorem~\ref{thm:firsteffective} gives
$$
 0=\Phi'(0)\le e^{-\frac{n-k-p+1}{p-1}t}\Phi'(t)\le0.
$$
Choose a sufficiently large regular value $t$ such that the
boundary collar lies in $\{0\le\psi<t\}$. Choose the first cutoff
parameter sufficiently small that the cutoff equals one on
$\Sigma_t$ and throughout a smaller closed boundary collar.
The integrated
divergence in the proof of Theorem~\ref{thm:firsteffective} is then
zero: its boundary flux is
$$
 (p-1)\bigl[e^{-\frac{n-k-p+1}{p-1}t}\Phi'(t)-\Phi'(0)\bigr]=0.
$$
The integrand is a sum of the nonnegative term containing
$\Div_gX$ and the nonnegative cutoff term. On the smaller collar,
the cutoff is one and its derivative is zero. If $\Div_gX$ were
positive at one point there, continuity would give a positive integral
on a small open set, contradicting the zero total integral.
Thus $\Div_gX=0$ on the collar and, by continuity, on the boundary.

Every summand in \eqref{eq:divX} is nonnegative. On
$\partial\Omega$, the square and the scalar curvature expression
have strictly positive coefficients. Since $\Div_gX=0$ there,
both expressions vanish. Using $\HH=h_\delta$ on the boundary
gives
$$
 \begin{aligned}
 (n-k)|\nabla\psi|-k\frac{\sigma_k(h_\delta)}{\sigma_{k-1}(h_\delta)}&=0,\\
 \frac{k(n-k-1)}{n-k}
 \frac{\sigma_k(h_\delta)^2}{\sigma_{k-1}(h_\delta)}
 -(k+1)\sigma_{k+1}(h_\delta)&=0.
 \end{aligned}
$$
On the boundary, strict $k$-convexity gives
$\sigma_{k-1}>0$ and $\sigma_k>0$. The second equality therefore
forces $\sigma_{k+1}>0$, so the curvature vector belongs to
 $\Gamma_{k+1}$. 
The second identity gives equality in the normalized Newton
inequality on $\Gamma_{k+1}$, so all principal curvatures coincide:
$h_\delta=\kappa I$. Since $\sigma_1>0$, their common value
$\kappa$ is positive. Consequently,
$$
 \frac{\sigma_k(h_\delta)}{\sigma_{k-1}(h_\delta)}
 =\frac{\binom{n-1}{k}\kappa^k}{\binom{n-1}{k-1}\kappa^{k-1}}
 =\frac{n-k}{k}\kappa.
$$
The first equality then says $|\nabla\psi|=\kappa$.

Strict star-shapedness makes the boundary a radial graph over
$\Sphere^{n-1}$, so it is connected. Let $\gamma$ be its induced
Euclidean metric.
Codazzi and $(h_\delta)_{ij}=\kappa\gamma_{ij}$ give
$$
 (\nabla_i\kappa)\gamma_{j\ell}
 =(\nabla_j\kappa)\gamma_{i\ell}.
$$
Contract with $\gamma^{j\ell}$. Since the boundary dimension is
$n-1$, the two sides become $(n-1)\nabla_i\kappa$ and
$\nabla_i\kappa$. Thus $(n-2)\nabla_i\kappa=0$, and
$\kappa$ is constant because $n\ge3$.
For any tangent vector $V$, our shape convention gives
$D_V\nu_\delta=h_\delta V=\kappa V$, hence
$$
 D_V\left(x-\frac{\nu_\delta}{\kappa}\right)
 =V-\kappa^{-1}\kappa V=0.
$$
The center $x_0=x-\nu_\delta/\kappa$ is therefore constant on
the connected boundary, and $|x-x_0|=1/\kappa$. The boundary
is a closed embedded hypersurface of the same dimension as this
sphere. Its inclusion into the sphere is locally open, and its
image is closed by compactness. The connected sphere is therefore
the entire boundary, so $\Omega$ is its bounded ball.
Only the Newton equality case on $\Gamma_{k+1}$ at the boundary
has been used. Convexity of the interior level hypersurfaces is
not required.
Together with the ball case proved above, this establishes the stated
equivalence.
\end{proof}

The proof separates the roles of the analytic hypotheses.
The cutoff comparisons between regular levels use
\textup{(P1)--(P3)}, while \textup{(P5)} supplies the asymptotic
information needed for the endpoint inequalities.
Assumption \textup{(P4)} gives weak flux conservation, which,
together with \textup{(P5)}, yields
$\CC(\Omega;u)=\mathsf a^{p-1}$.
The relation between ellipticity and \textup{(P3)} is discussed
in Appendix~\ref{app:ellipticity}.

\section{The limit as \texorpdfstring{\(p\downarrow1\)}{p tends to 1}}
\label{sec:limit}
We use the construction in \cite{ExteriorConstruction} and the
fixed-parameter inequality to derive the geometric Minkowski inequality.

\begin{proposition}\label{ass:comparison}
Let $n\ge3$, $1\le k\le n-2$, and let $\Omega$ be a smooth bounded
strictly star-shaped domain with strictly $k$-convex boundary.
There exists $\varepsilon_0=\varepsilon_0(\Omega,n,k)\in(0,n-k-1)$
such that, for every $1<p<1+\varepsilon_0$, there is a positive
exterior potential $u_p$ satisfying
Assumptions~\ref{ass:potential} and~\ref{ass:asymptotic}.
Moreover, there exists $r_*=r_*(\Omega)>0$ such that, for every
$x\in\partial\Omega$, the ball $B_{r_*}(x-r_*\nu(x))$ lies in
$\Omega$ and its potential satisfies
\[
 v_{p,x}(y)
 =\left(\frac{r_*}{|y-x+r_*\nu(x)|}\right)^{
                         \frac{n-k-p+1}{p-1}}
 \le u_p(y)\qquad(y\in E).
\]
\end{proposition}

\begin{proof}
 Take $\varepsilon_0$ from \cite[Theorem~1.1]{ExteriorConstruction}. For every
 $1<p<1+\varepsilon_0<n-k$, that theorem supplies a potential $u_p$
 satisfying \textup{(P1)--(P5)}.

 The compact smooth boundary has a uniform interior sphere radius
 $r_*>0$. This radius depends only on the fixed domain. For each
 $x\in\partial\Omega$, apply the interior-ball comparison proved in
 \cite[Proposition~2.2]{ExteriorConstruction} to $B_{r_*}(x-r_*\nu(x))$ and the
 constructed potential $u_p$. The comparison allows tangency
 at $x$ and gives the stated inequality throughout the exterior.
 \end{proof}

The far-field estimates in \textup{(P5)} are used at each fixed $p$
and need not be uniform as $p\downarrow1$.

\begin{remark*}
A proper regularization construction on general smooth domains,
together with comparison for the selected potentials, will be
treated in a separate paper. For arbitrary potentials satisfying
\textup{(P1)--(P5)}, interior-ball comparison must be established
separately, since these hypotheses specify neither a viscosity
branch nor a selection procedure.
\end{remark*}

\begin{proposition}%
\label{prop:limit}
Every smooth bounded strictly star-shaped domain $\Omega$ with
strictly $k$-convex boundary satisfies
\begin{equation}\label{eq:AFtarget}
 A_k(\Omega)\ge A_{k-1}(\Omega)^{(n-k-1)/(n-k)}.
\end{equation}
Equality holds if and only if $\Omega$ is a ball.
\end{proposition}

\begin{proof}
Let $u_p$ be the potentials supplied by
Proposition~\ref{ass:comparison}, and set
\[
 \varepsilon=p-1,\qquad
 q_p=\frac{\varepsilon|Du_p|}{n-k-p+1}
 \quad\text{on }\partial\Omega.
\]
We first establish positive upper and lower bounds for $q_p$
independent of sufficiently small $\varepsilon$.

Write
$\partial\Omega=\{r_\Omega(\theta)\theta:\theta\in\Sphere^{n-1}\}$
and define the radial gauge
\[
 G(r\theta)=r/r_\Omega(\theta).
\]
The function $G$ is smooth away from the origin, homogeneous of
degree one, and has nonzero gradient. Its level hypersurfaces are
homothetic copies of $\partial\Omega$.

For $\xi\ne0$, set
\[
 P_p(\xi)=I+\left(\sqrt{\frac{\varepsilon}{k}}-1\right)
                 \frac{\xi\otimes\xi}{|\xi|^2},
 \qquad
 M_p[v]=-P_p(Dv)D^2vP_p(Dv).
\]
By the weighted-Hessian identity in
Section~\ref{sec:weightedhessian} and \textup{(P2)},
$\sigma_k(M_p[u_p])=0$ at noncritical points.
On $G=1$, put $e=DG/|DG|$ and
\[
 D^2G|_{e^\perp}=|DG|h_\delta.
\]
For $1\le j\le k$, define
\[
 a_j=\sigma_j\bigl(D^2G|_{e^\perp}\bigr),\qquad
 b_j=|DG|^2\sigma_{j-1}\bigl(D^2G|_{e^\perp}\bigr),\qquad
 d_j=\sigma_j(D^2G)-a_j.
\]
Strict $k$-convexity and compactness give positive lower bounds
for $a_j$ and $b_j$.

For a fixed small $c_*>0$, let
\[
 m_p=\frac{c_*}{\varepsilon},\qquad w_p=G^{-m_p}.
\]
Then
\[
 M_p[w_p]=m_pG^{-m_p-1}Q_p,\qquad
 Q_p=P_p(DG)D^2G P_p(DG)
     -\frac{\varepsilon}{k}(m_p+1)G^{-1}DG\otimes DG.
\]
In an orthonormal frame with last vector $e$, conjugation by
$P_p(DG)$ multiplies each principal minor involving the last
index by $\varepsilon/k$, while the rank-one term changes only
the last diagonal entry. Thus, on $G=1$,
\[
 \sigma_j(Q_p)
 =a_j+\frac{\varepsilon}{k}d_j
      -\frac{\varepsilon}{k}(m_p+1)b_j
 =a_j-\frac{c_*}{k}b_j+\frac{\varepsilon}{k}(d_j-b_j).
\]
Choose $c_*$ so that $(c_*/k)b_j\le a_j/4$ for every
$j\le k$ and every boundary point. For sufficiently small
$\varepsilon$, the last term has absolute value at most $a_j/4$,
and hence $\sigma_j(Q_p)\ge a_j/2>0$.
Since $Q_p$ is homogeneous of degree $-1$, these inequalities
hold throughout $E$. Consequently,
\[
 M_p[w_p]\in\Gamma_k,\qquad Dw_p\ne0,
 \qquad w_p=1\text{ on }\partial\Omega,\qquad
 w_p\longrightarrow0\text{ at infinity}.
\]

If $u_p-w_p$ were positive somewhere, it would attain a positive
maximum at a finite interior point. There
$Du_p=Dw_p\ne0$, so \textup{(P2)} applies, and Hessian ordering gives
\[
 M_p[u_p]\ge M_p[w_p]\in\Gamma_k.
\]
Since $\Gamma_k$ is preserved by addition of positive semidefinite
matrices, this contradicts $\sigma_k(M_p[u_p])=0$.
Thus $u_p\le w_p$.
Differentiating at the boundary, where
$DG=|DG|\nu_\delta$ and $\partial_{\nu_\delta} u_p=-|Du_p|$, yields
\[
 |Du_p|\ge m_p|DG|,
 \qquad
 q_p\ge
 \frac{c_*\min_{\partial\Omega}|DG|}{n-k}
 >0.
\]

For the upper bound, Proposition~\ref{ass:comparison} gives
$v_{p,x}\le u_p$, with equality at $x\in\partial\Omega$.
Differentiation along the outward normal gives
\[
 -|Du_p(x)|=\partial_{\nu_\delta} u_p(x)
 \ge\partial_{\nu_\delta} v_{p,x}(x)
 =-\frac{n-k-p+1}{\varepsilon r_*}.
\]
Therefore
\[
 0<\frac{c_*\min_{\partial\Omega}|DG|}{n-k}\le q_p\le r_*^{-1}<\infty
 \qquad\text{on }\partial\Omega,
\]
uniformly for sufficiently small $\varepsilon$.

These bounds imply $|\log q_p|\le C$ and hence
\[
 \sup_{\partial\Omega}|q_p^{p-1}-1|
 \le e^{C(p-1)}-1\longrightarrow0.
\]
By the definition of the normalized flux,
\[
 \CC(\Omega;u_p)
 =\avgint_{\partial\Omega}q_p^{p-1}H_{k-1}\dd S_\delta
 \longrightarrow A_{k-1}(\Omega).
\]
Likewise, $H_k/H_{k-1}$ has positive upper and lower bounds on
the fixed boundary, so
\[
 \frac{H_k^p}{H_{k-1}^{p-1}}
 =H_k\left(\frac{H_k}{H_{k-1}}\right)^{p-1}
 \longrightarrow H_k
 \quad\text{uniformly}.
\]
Since $A_{k-1}>0$, it follows that
\[
 \CC(\Omega;u_p)^{\frac{n-k-p}{n-k-p+1}}
 \longrightarrow A_{k-1}^{\frac{n-k-1}{n-k}},
 \qquad
 \avgint_{\partial\Omega}
 \frac{H_k^p}{H_{k-1}^{p-1}}\dd S_\delta
 \longrightarrow A_k.
\]
Passing to the limit in \eqref{eq:fixedinequality} proves
\eqref{eq:AFtarget}.
The far-field hypothesis \textup{(P5)} is used at each fixed $p$;
its remainder need not be uniform.
Once the potential family and its comparison property are
available, this limit passage requires neither convergence of
the logarithmic potentials nor uniform regularity estimates in $p$.

Suppose equality holds in \eqref{eq:AFtarget}, and set
\[
 \mathcal F(\Omega)
 =A_k-A_{k-1}^{(n-k-1)/(n-k)},\qquad
 \Lambda=A_{k-1}^{-1/(n-k)}.
\]
Strict star-shapedness and strict $k$-convexity persist under
sufficiently small $C^2$ normal variations
$x_t=x+t f(x)\nu_\delta(x)$, with either sign of $t$.
The inequality just proved applies to each resulting domain
$\Omega_t$, giving
$\mathcal F(\Omega_t)\ge0=\mathcal F(\Omega)$.
The parameter threshold may depend on $\Omega_t$.
Thus $\delta\mathcal F=0$ for every smooth normal speed $f$.

With the convention $h_\delta=D\nu_\delta$, the first-variation formula is
\[
 \delta A_j
 =\frac{n-1-j}{\omega_{n-1}}
   \int_{\partial\Omega}H_{j+1}f\dd S_\delta.
\]
For $j\ge1$, this follows from
$\dot h_\delta=-\nabla_\Sigma^2f-fh_\delta^2$, the area variation,
the divergence-free Newton tensor, and the identity
\[
 \tr\bigl(T_{j-1}(h_\delta)h_\delta^2\bigr)=\sigma_1\sigma_j-(j+1)\sigma_{j+1}.
\]
The case $j=0$ is the area variation itself.
Consequently,
\[
 0=\delta\mathcal F
 =\frac{n-k-1}{\omega_{n-1}}
   \int_{\partial\Omega}
   (H_{k+1}-\Lambda H_k)f\dd S_\delta.
\]
Since $f$ is arbitrary and $n-k-1>0$,
\[
 H_{k+1}=\Lambda H_k.
\]

As $\Lambda>0$ and $H_k>0$, the curvature vector belongs to
$\Gamma_{k+1}$. Newton's inequality gives
\[
 H_k^2\ge H_{k-1}H_{k+1}
 =\Lambda H_{k-1}H_k,
 \qquad H_k\ge\Lambda H_{k-1}.
\]
On the other hand,
\[
 \avgint_{\partial\Omega}
 (H_k-\Lambda H_{k-1})\dd S_\delta
 =\mathcal F(\Omega)=0.
\]
Hence $H_k=\Lambda H_{k-1}$ everywhere, and equality holds
in Newton's inequality. All principal curvatures therefore
coincide. Writing $h_\delta=\kappa I$, the relation
$H_{k+1}=\Lambda H_k$ gives $\kappa=\Lambda>0$.
Thus
\[
 D_V(x-\Lambda^{-1}\nu_\delta)=0
\]
for every tangent vector $V$. Since $\partial\Omega$ is connected
and closed, it is a sphere.
Conversely, $A_j(B_R)=R^{n-1-j}$ gives equality for every ball.
\end{proof}

\clearpage
\appendix
\section{Ellipticity and level-curvature positivity}\label{app:ellipticity}

We examine the linearization of the operator and its relation to
the level-curvature hypothesis \textup{(P3)}. All derivatives,
inner products, and norms in this appendix are Euclidean.
The outward flux \(F_z\) and the matrix \(P\) are those defined in
Section~\ref{sec:weightedhessian}, so that
\(\Div F_z=(-1)^k\mathcal L_{k,p}u\).
All identities involving negative powers of \(|Du|\) are restricted
to the noncritical set. The calculation for radial potentials
also records the dependence of the ellipticity ratio on \(p\).
\subsection{The principal matrix and level-curvature positivity}

For the linearization we use the equivalent equation
\(-\Div F_z=0\), so that its principal matrix has the positive
elliptic sign on the strictly elliptic branch.

\begin{proposition}[Principal linearization]\label{prop:principal}
Holding \(Du\) fixed, the derivative of \(-\Div F_z\) with respect
to \(D^2u\) is represented by
\begin{equation}\label{eq:principal}
 \mathsf A
 =k|Du|^{p-k-1}P T_{k-1}(-PD^2uP)P.
\end{equation}
Consequently, strict ellipticity at a noncritical point is equivalent
to \(T_{k-1}(-PD^2uP)>0\), and
\begin{equation}\label{eq:normalprincipal}
 \mathsf A(\nu,\nu)
 =(p-1)|Du|^{p-2}\sigma_{k-1}(h_\delta).
\end{equation}
In particular, strict ellipticity supplies the pointwise curvature
positivity required in Assumption~\ref{ass:potential}\textup{(P3)}.
\end{proposition}
\begin{proof}
\medskip
\noindent\emph{Step 1. The principal linearization.}\quad
The quantities \(|Du|,\nu,P\) are held fixed. Their derivatives
contribute only to the first-order part of the linearization.
For a symmetric Hessian variation \(H\), differentiating
\(-k|Du|^{p-k-1}\sigma_k(-PD^2uP)\) gives
\[
 \begin{aligned}
 D_{D^2u}(-\Div F_z)[H]
 &=k|Du|^{p-k-1}
     \tr\!\left(T_{k-1}(-PD^2uP)PHP\right)\\
 &=\tr(\mathsf A H).
 \end{aligned}
\]
The two minus signs cancel: one comes from \(-\Div F_z\), the other
from the Hessian argument \(-PD^2uP\).
The cyclic trace identity places both copies of \(P\) around the
Newton tensor. For every vector \(\eta\),
\[
 \mathsf A(\eta,\eta)
 =k|Du|^{p-k-1}T_{k-1}(-PD^2uP)(P\eta,P\eta).
\]
The scalar factor is positive, and \(P\) is invertible for \(p>1\).
Thus \(\mathsf A>0\) is equivalent to
\(T_{k-1}(-PD^2uP)>0\).

\medskip
\noindent\emph{Step 2. The normal coefficient.}\quad
With the outward normal \(\nu=-Du/|Du|\), tangent vectors \(V,W\)
satisfy
\[
 -D^2u(V,W)=|Du|h_\delta(V,W).
\]
Since \(P\) fixes tangent vectors, the tangential block of
\(-PD^2uP\) is \(|Du|h_\delta\). The normal entry of its Newton
tensor is the sum of the principal \((k-1)\)-minors avoiding that
index, so
\[
 \bigl(T_{k-1}(-PD^2uP)\bigr)_{\nu\nu}
 =|Du|^{k-1}\sigma_{k-1}(h_\delta),\qquad
 P\nu=\sqrt{\frac{p-1}{k}}\,\nu.
\]
Therefore
\[
 \begin{aligned}
 \mathsf A(\nu,\nu)
 &=k|Du|^{p-k-1}\frac{p-1}{k}
       |Du|^{k-1}\sigma_{k-1}(h_\delta)\\
 &=(p-1)|Du|^{p-2}\sigma_{k-1}(h_\delta).
 \end{aligned}
\]
If \(\mathsf A>0\), division by \((p-1)|Du|^{p-2}>0\)
gives \textup{(P3)} at this point. The converse does not follow from this
single matrix entry.
\end{proof}

For a \(C^2\) solution, pointwise strict positivity of \(\mathsf A\)
gives uniform ellipticity on each compact subset of \(\{Du\ne0\}\):
its smallest eigenvalue has a positive minimum there, and its
largest eigenvalue has a finite maximum. The constants may depend
on the compact subset, the selected solution, and \(p\); they are
not asserted to be uniform as \(p\downarrow1\).

If the solution is already \(C^{2,\alpha}\), \(0<\alpha<1\),
the usual local Schauder bootstrap applies on these compact
noncritical sets. Differentiating \(-\Div F_z=0\) gives linear
equations for the first derivatives of \(u\), with principal matrix
\(\mathsf A\) and H\"older continuous coefficients. Difference
quotients justify the differentiation, and interior Schauder
estimates give \(C^{3,\alpha}\) regularity and then further
bootstrapping. This argument supplies neither the initial
\(C^{2,\alpha}\) estimate nor boundary regularity.

The equation also determines the normal Hessian whenever
\(\sigma_{k-1}(h_\delta)>0\). Write \(D_\Sigma\) for the Euclidean
tangential gradient along a level set. Since
\(-D^2u(V,\nu)=D_V|Du|\), the weighted Hessian in a
level-adapted orthonormal frame is
\[
 -PD^2uP=
 \begin{pmatrix}
 |Du|h_\delta&\sqrt{(p-1)/k}\,D_\Sigma|Du|\\
 \sqrt{(p-1)/k}\,(D_\Sigma|Du|)^{\mathsf T}
      &-\dfrac{p-1}{k}u_{\nu\nu}
 \end{pmatrix}.
\]
Expanding the principal minors and using homogeneity gives
\[
 \begin{aligned}
 \sigma_k(-PD^2uP)
 &=|Du|^k\sigma_k(h_\delta)
   -\frac{p-1}{k}u_{\nu\nu}|Du|^{k-1}\sigma_{k-1}(h_\delta)\\
 &\quad-\frac{p-1}{k}|Du|^{k-2}
      T_{k-2}(h_\delta)(D_\Sigma|Du|,D_\Sigma|Du|).
 \end{aligned}
\]
The equation says that this expression is zero. Solving for the
normal entry and dividing by
\(\frac{p-1}{k}|Du|^{k-1}\sigma_{k-1}(h_\delta)>0\) yields
\begin{equation}\label{eq:normalhessian}
 u_{\nu\nu}
 =\frac{k}{p-1}|Du|
       \frac{\sigma_k(h_\delta)}{\sigma_{k-1}(h_\delta)}
  -\frac{T_{k-2}(h_\delta)(D_\Sigma|Du|,D_\Sigma|Du|)}
        {|Du|\sigma_{k-1}(h_\delta)}.
\end{equation}
For \(k=1\), \(T_{-1}=0\) by convention, and the formula becomes
\((p-1)u_{\nu\nu}=|Du|\sigma_1(h_\delta)\), as expected for
the \(p\)-Laplace equation. The identity is a starting point for
Hessian estimates, but gives no quantitative control of its
denominators or the tangential quantities on the right.

\begin{remark}[A local obstruction to automatic strictness]
For \(k\ge2\), consider near the origin
\[
 u(x)=\frac12-x_n-\frac{x_n^2}{2}
                  -\frac12\sum_{i=1}^{k-2}x_i^2,
\]
where the sum is empty for \(k=2\). This function is positive and
has nonzero gradient in a small neighborhood. The matrix
\(-D^2u\) is positive semidefinite of rank \(k-1\), and
\(-PD^2uP\) has the same properties. Thus
\[
 \begin{gathered}
 -PD^2uP\in\overline{\Gamma_k},\\
 \sigma_{k-1}(-PD^2uP)>0,\qquad \sigma_k(-PD^2uP)=0,
 \end{gathered}
\]
and \(\Div F_z=0\) there. At the origin the tangential Hessian
has rank \(k-2\), hence \(\sigma_{k-1}(h_\delta)=0\).
Even noncritical smoothness and strict ambient
\(\sigma_{k-1}(-PD^2uP)>0\) therefore do not imply
Assumption~\ref{ass:potential}\textup{(P3)} or strict ellipticity.
This is a local example, not a counterexample to the exterior
problem with its prescribed boundary and infinity conditions.
\end{remark}

\subsection{Ellipticity in the spherical model}

For a ball of radius \(R\), the positive potential is
\[
 u(x)=\left(\frac Rr\right)^{\frac{n-k-p+1}{p-1}},
 \qquad r=|x|.
\]
Its radial derivatives satisfy
\[
 |Du|=-u_r=\frac{n-k-p+1}{p-1}\frac{u}{r}>0,\qquad
 u_{rr}=\frac{n-k}{p-1}\frac{|Du|}{r}.
\]
The tangential eigenvalues of \(-D^2u\) are \(|Du|/r\);
the normal eigenvalue is \(-u_{rr}\). Multiplication by \(P\)
on both sides scales only the normal entry by \((p-1)/k\), so
\[
 -PD^2uP
 =\frac{|Du|}{r}
    \operatorname{diag}\!\left(I_{n-1},-\frac{n-k}{k}\right).
\]
The normal eigenvalue of \(T_{k-1}(-PD^2uP)\) is
\[
 \binom{n-1}{k-1}\left(\frac{|Du|}{r}\right)^{k-1}.
\]
For a tangential direction, omitting that direction from the
principal minors instead gives
\[
 \begin{aligned}
 &\left(\frac{|Du|}{r}\right)^{k-1}
  \left[\binom{n-2}{k-1}
        -\frac{n-k}{k}\binom{n-2}{k-2}\right]\\
 &\qquad=\binom{n-1}{k-1}
       \left(\frac{|Du|}{r}\right)^{k-1}
       \frac{n-k}{k(n-1)}.
 \end{aligned}
\]
Here we used
\(\binom{n-2}{k-1}=\frac{n-k}{n-1}\binom{n-1}{k-1}\)
and
\(\binom{n-2}{k-2}=\frac{k-1}{n-1}\binom{n-1}{k-1}\).
For \(k=1\), a binomial coefficient with lower index \(-1\)
is zero.

In \eqref{eq:principal}, the two factors of \(P\) contribute
an additional factor \((p-1)/k\) in the normal direction and
no change in tangential directions. Since
\[
 |Du|^{p-k-1}
       \left(\frac{|Du|}{r}\right)^{k-1}
 =|Du|^{p-2}r^{1-k},
\]
the principal matrix is
\begin{equation}\label{eq:balllinearization}
 \mathsf A=\binom{n-1}{k-1}|Du|^{p-2}r^{1-k}
 \left[\frac{n-k}{n-1}(I-\nu\otimes\nu)
                    +(p-1)\nu\otimes\nu\right].
\end{equation}
The normal-to-tangential coefficient ratio is
\[
 \frac{(p-1)(n-1)}{n-k}\longrightarrow0
 \qquad\text{as }p\downarrow1.
\]
This ratio is unchanged by scalar normalization of the equation.
Thus the ellipticity ratio is not uniform as \(p\downarrow1\),
even for balls. The boundary-flux limit in
Section~\ref{sec:limit} does not require uniform interior
Schauder estimates for the parameter family.

\section{The first variation of the energy}\label{app:energy-detailed}

For a function \(w\) on an open set \(U\subset\RR^n\), write
\begin{equation}\label{eq:energy-definition-detailed}
 I_{k,p}[w]=\int_U |Dw|^{p-k-1}
             \ip{T_{k-1}(D^2w)Dw}{Dw}\dd x
\end{equation}
whenever the integral exists.

\begin{proposition}[First variation]
Let \(w\in C^3(U)\). For every smooth variation \(\varphi\) compactly
supported in \(\{|Dw|>0\}\), the local first variation is
\begin{equation}\label{eq:energy-variation-detailed}
 \delta I_{k,p}[w](\varphi)
 =p\int_{\{|Dw|>0\}} |Dw|^{p-k-1}\ip{T_{k-1}(D^2w)Dw}{D\varphi}\dd x
 =-p\int_{\{|Dw|>0\}}\varphi\,\LL w\dd x .
\end{equation}
\end{proposition}

\begin{proof}
For \(k=1\), differentiate \(\int_U|Dw|^p\dd x\). For \(k\ge2\), put
\(b(s)=\int_1^s t^{p-k}\dd t\) for \(s>0\), so \(b(s)=\log s\) if \(p=k-1\).
Since \(T_{k-2}(D^2w)\) is divergence-free and
\(\tr(T_{k-2}(D^2w)D^2w)=(k-1)\sigma_{k-1}(D^2w)\), the Newton
recurrence and \(s^{p-k+1}=1+(p-k+1)b(s)\) give
\[
\begin{aligned}
 |Dw|^{p-k-1}\ip{T_{k-1}(D^2w)Dw}{Dw}
 ={}&\bigl(1+p\,b(|Dw|)\bigr)\sigma_{k-1}(D^2w)\\
 &-\Div\!\bigl(b(|Dw|)\,T_{k-2}(D^2w)Dw\bigr).
\end{aligned}
\]
The divergence term has zero first variation for this compactly supported
perturbation. In the remaining term the variation of \(\sigma_{k-1}\) is
\(\tr(T_{k-2}(D^2w)D^2\varphi)\). Integrating this term once by parts gives
\[
 \delta I_{k,p}[w](\varphi)
 =p\int_{\{|Dw|>0\}}\ip{|Dw|^{p-k-1}\sigma_{k-1}(D^2w)Dw
     -T_{k-2}(D^2w)D[b(|Dw|)]}{D\varphi}\dd x .
\]
Now \(D[b(|Dw|)]=|Dw|^{p-k-1}D^2w\,Dw\). The Newton recurrence gives
the first equality in \eqref{eq:energy-variation-detailed}, and a final
integration by parts gives the second.
\end{proof}

\begin{proposition}[Exterior energy and boundary flux]
If \(u\) satisfies \textup{(P1)--(P5)}, then
\begin{equation}\label{eq:energy-boundary-identity-detailed}
 (-1)^{k-1}I_{k,p}[u]
 =\int_E |Du|^p\sigma_{k-1}(h_\delta)\dd x
 =\int_{\partial\Omega}|Du|^{p-1}\sigma_{k-1}(h_\delta)\dd S_\delta,
\end{equation}
where the volume density is extended by zero on \(\{Du=0\}\).
\end{proposition}

\begin{proof}
At a regular point \(D^2u|_{\nu^\perp}=-|Du|h_\delta\), with
\(\nu=-Du/|Du|\). Taking the normal component of the Newton tensor gives
\begin{equation}\label{eq:energy-density-sign}
 (-1)^{k-1}|Du|^{p-k-1}\ip{T_{k-1}(D^2u)Du}{Du}
 =|Du|^p\sigma_{k-1}(h_\delta).
\end{equation}
The extended density thus equals the scalar product of \(-Du\) with
the divergence-free Newton field of \textup{(P4)} almost everywhere.
For a regular value \(0<t<1\), apply the divergence theorem to
\(t-u\) times this field on \(\{t<u<1\}\). The distributional
product rule applies since \(u\in C^1\), and the field is classical
in both boundary collars.
The multiplier vanishes on \(\{u=t\}\). On \(\partial\Omega\) it
equals \(t-1\) and the integration domain has outward normal \(-\nu\).
Hence \eqref{eq:flatflux} gives
\begin{equation}\label{eq:energy-annulus-ibp}
 \int_{\{t<u<1\}}|Du|^p\sigma_{k-1}(h_\delta)\dd x
 =(1-t)\int_{\partial\Omega}|Du|^{p-1}\sigma_{k-1}(h_\delta)\dd S_\delta .
\end{equation}
By \textup{(P5)} all sufficiently small \(t>0\) are regular.
The density is nonnegative by \textup{(P3)}, so monotone convergence
as \(t\downarrow0\) proves \eqref{eq:energy-boundary-identity-detailed}.
\end{proof}

\bibliographystyle{plain}
\bibliography{references/paper1/reference}
\end{document}